\documentclass[12pt]{amsart}
\usepackage{amsmath, amsthm, amscd, amsfonts, amssymb, graphicx, color, alltt}
\usepackage[bookmarksnumbered, plainpages]{hyperref}
\usepackage{alltt}
\usepackage{enumerate}
\usepackage{subfig}
\newtheorem{thm}{Theorem}[section]

\newtheorem{lem}[thm]{Lemma}

\newtheorem{conj}[thm]{Conjecture}

\numberwithin{equation}{section}
\usepackage[utf8]{inputenc}
\usepackage{pgf,tikz}
\usetikzlibrary{arrows}

\title{\bf Counting Centralizers of a Finite Group with an Application in Constructing the Commuting Conjugacy Class Graph}

\author{\bf  A. R. Ashrafi and M. A. Salahshour$^\star$}

\thanks{$^\star$Corresponding author (Email: salahshour@iausk.ac.ir)}

\address{\textbf{Ali Reza Ashrafi:} Department of Pure Mathematics, Faculty of Mathematical Sciences, University of Kashan, Kashan 87317$-$53153, I. R. Iran}

\address{\textbf{Mohammad Ali Salahshour:} Department of Mathematics, Savadkooh Branch, Islamic Azad University, Savadkooh, I. R. Iran}

\date{}
\begin{document}

\maketitle

\begin{abstract}
The set of all centralizers of elements in a finite group $G$ is denoted by $Cent(G)$ and $G$ is called $n-$centralizer if $|Cent(G)| = n$. In this paper, the structure of centralizers in a non-abelian finite group $G$ with this property that $\frac{G}{Z(G)} \cong Z_{p^2} \rtimes Z_{p^2}$ is obtained. As a consequence, it is proved that such a group has exactly $[(p+1)^2+1]$ element centralizers and  the structure of the commuting conjugacy class graph of $G$ is completely determined.

\vskip 3mm


\vskip 3mm

\end{abstract}

\bigskip

\section{Introduction}
Throughout this paper all groups are assumed to be finite and $C_G(x)$ denotes the  centralizer of an element $x$ in $G$. If $G$ is a group containing two subgroups $H$ and $N$ such that $N\lhd G$, $H \cap N = 1$ and $G=HN$ then $G$ is said to be the semi-direct product of $H$ by $N$ and we write $G = N \rtimes H$. The group $G$ is called capable if there exists another group $H$ such that $G \cong \frac{H}{Z(H)}$. Our other notations are standard and taken mainly from \cite{5} and our calculations are done with the aid of GAP \cite{7}.

The center of $G$ is denoted by $Z(G)$ and  $Cent(G)$ = $\{ C_G(x) \mid x \in G\}$. The group $G$ is called $n-$centralizer if $n = |Cent(G)|$ \cite{1}.  The study of finite groups with respect to the number of distinct element centralizers was started by Belcastro and Sherman \cite{2}. It is clear that  $|Cent(G)| = 1$  if and only if $G$ is abelian and there is no finite group with exactly two or three element centralizers. They proved that if $\frac{G}{Z(G)} \cong \mathbb{Z}_p \times \mathbb{Z}_p$  then $|Cent(G)| = p + 2$ \cite[Theorem 5]{2}.

Let  $H$ be a graph with vertex set $\{ 1, 2, \cdots, k\}$ and $G_i$'s, $1 \leq i \leq k$, are disjoint
graphs of order $n_i$. Following Sabidussi \cite{6}, the graph $H[G_1, G_2,
\cdots, G_k]$ is formed by taking the graphs $G_1$, $G_2$,
$\cdots,$ $G_k$ and connect a vertex of $G_i$ to another vertex
in $G_j$ whenever $i$ is adjacent to $j$ in $H$. The graph
$H\left[ G_1, G_2, \cdots, G_k\right]$ is called the $H-$join of
the graphs $G_1, \cdots, G_k$.

The notion of \textbf{commuting conjugacy class graph} of a non-abelian group $G$, $\Gamma(G)$, was introduced by Mohammadian et al.  \cite{4}.   This is a simple graph with non-central conjugacy classes of $G$ as its vertex  set and two distinct vertices $A$ and $B$ are adjacent if and only if there are $x \in A$ and $y \in B$ such that $xy = yx$. In such a case, we also say that two conjugacy classes commute to each other. The authors of the mentioned paper obtained some interesting properties of this graph among them a classification of  triangle-free commuting conjugacy class graph of a finite groups is given.

A \textbf{CA-group} is a group in which every noncentral element has an abelian centralizer. The present authors \cite{65}, obtained the structure of the commuting conjugacy class graph of finite CA-groups. They proved that this graph is a union of some complete graphs. As a consequence of this work, the commuting conjugacy class graph of dihedral, semi-dihedral, dicyclic and three other families of meta-cyclic groups were constructed.

The following two theorems are the main results of this paper:

\begin{thm}\label{t1}
Suppose $p$ is a prime number and $G$ is a group with center $Z$ such that $\frac{G}{Z} \cong \mathbb{Z}_{p^2} \rtimes \mathbb{Z}_{p^2}$.  Then $G$ has exactly $[(p+1)^2+1]$ element centralizers.
\end{thm}

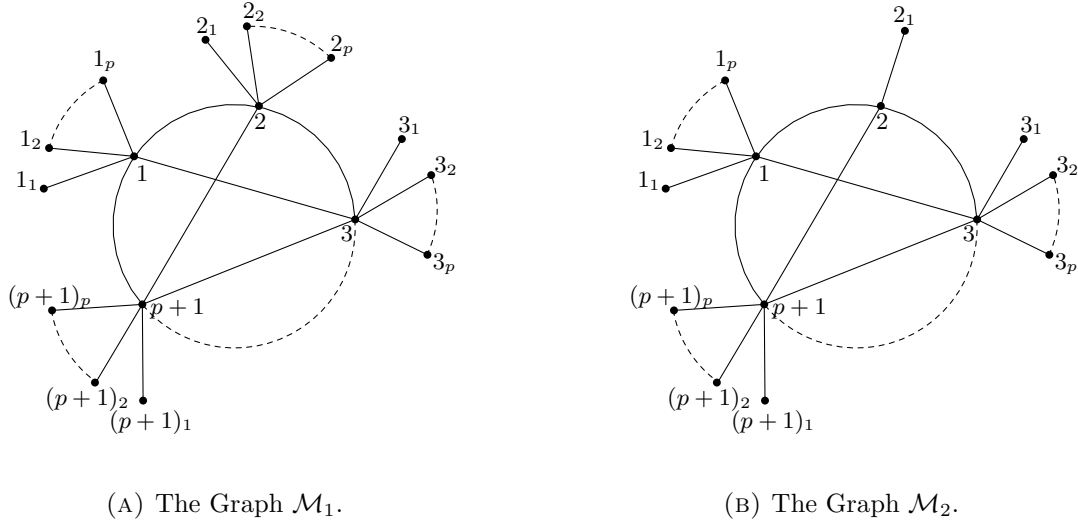
\begin{figure}[h]
\centering
\subfloat[The Graph $\mathcal{M}_1$.]{
\begin{picture}(200,200)(75,75)
\definecolor{qqqqff}{rgb}{0,0,0}
\begin{tikzpicture}[line cap=round,line join=round,>=triangle 45,x=1cm,y=1cm]
\clip(-5.66,-4.54) rectangle (15.68,4.94);
\draw [shift={(0.55,1.41)}] plot[domain=0.05:3.86,variable=\t]({1*1.61*cos(\t r)+0*1.61*sin(\t r)},{0*1.61*cos(\t r)+1*1.61*sin(\t r)});
\draw [shift={(0.56,1.38)},dash pattern=on 2pt off 2pt] plot[domain=-2.44:0.07,variable=\t]({1*1.6*cos(\t r)+0*1.6*sin(\t r)},{0*1.6*cos(\t r)+1*1.6*sin(\t r)});
\draw [shift={(0.81,2.6)},dash pattern=on 2pt off 2pt] plot[domain=0.79:1.63,variable=\t]({1*1.46*cos(\t r)+0*1.46*sin(\t r)},{0*1.46*cos(\t r)+1*1.46*sin(\t r)});
\draw [shift={(1.86,1.62)},dash pattern=on 2pt off 2pt] plot[domain=-0.45:0.33,variable=\t]({1*1.39*cos(\t r)+0*1.39*sin(\t r)},{0*1.39*cos(\t r)+1*1.39*sin(\t r)});
\draw [shift={(-0.42,0.53)},dash pattern=on 2pt off 2pt] plot[domain=3.31:4.09,variable=\t]({1*1.48*cos(\t r)+0*1.48*sin(\t r)},{0*1.48*cos(\t r)+1*1.48*sin(\t r)});
\draw [shift={(-0.51,2.06)},dash pattern=on 2pt off 2pt] plot[domain=2.06:2.88,variable=\t]({1*1.45*cos(\t r)+0*1.45*sin(\t r)},{0*1.45*cos(\t r)+1*1.45*sin(\t r)});
\draw (-0.78,2.33)-- (-1.98,1.9);
\draw (-1.91,2.44)-- (-0.78,2.33);
\draw (-1.19,3.34)-- (-0.78,2.33);
\draw (0.17,3.88)-- (0.88,3);
\draw (0.72,4.06)-- (0.88,3);
\draw (1.84,3.64)-- (0.88,3);
\draw (-0.66,-0.92)-- (-0.67,0.36);
\draw (-1.29,-0.68)-- (-0.67,0.36);
\draw (-1.87,0.28)-- (-0.67,0.36);
\draw (0.88,3)-- (-0.67,0.36);
\draw (3.12,1.02)-- (2.16,1.49);
\draw (3.17,2.08)-- (2.16,1.49);
\draw (2.78,2.56)-- (2.16,1.49);
\draw (-0.78,2.33)-- (2.16,1.49);
\draw (2.16,1.49)-- (-0.67,0.36);
\begin{scriptsize}
\fill [color=qqqqff] (-0.78,2.33) circle (1.5pt);
\draw[color=qqqqff] (-0.67,2.1) node {$1$};
\fill [color=qqqqff] (-1.19,3.34) circle (1.5pt);
\draw[color=qqqqff] (-1.17,3.6) node {$1_p$};
\fill [color=qqqqff] (-1.91,2.44) circle (1.5pt);
\draw[color=qqqqff] (-2.15,2.55) node {$1_2$};
\fill [color=qqqqff] (-1.98,1.9) circle (1.5pt);
\draw[color=qqqqff] (-2.2,2.0) node {$1_1$};
\fill [color=qqqqff] (0.88,3) circle (1.5pt);
\draw[color=qqqqff] (0.9,2.8) node {$2$};
\fill [color=qqqqff] (1.84,3.64) circle (1.5pt);
\draw[color=qqqqff] (1.99,3.8) node {$2_p$};
\fill [color=qqqqff] (0.72,4.06) circle (1.5pt);
\draw[color=qqqqff] (0.82,4.25) node {$2_2$};
\fill [color=qqqqff] (0.17,3.88) circle (1.5pt);
\draw[color=qqqqff] (0.2,4.06) node {$2_1$};
\fill [color=qqqqff] (2.16,1.49) circle (1.5pt);
\draw[color=qqqqff] (2.05,1.3) node {$3$};
\fill [color=qqqqff] (3.12,1.02) circle (1.5pt);
\draw[color=qqqqff] (3.35,.9) node {$3_p$};
\fill [color=qqqqff] (3.17,2.08) circle (1.5pt);
\draw[color=qqqqff] (3.37,2.2) node {$3_2$};
\fill [color=qqqqff] (2.78,2.56) circle (1.5pt);
\draw[color=qqqqff] (2.89,2.74) node {$3_1$};
\fill [color=qqqqff] (-0.67,0.36) circle (1.5pt);
\draw[color=qqqqff] (-0.2,0.3) node {$p+1$};
\fill [color=qqqqff] (-1.87,0.28) circle (1.5pt);
\draw[color=qqqqff] (-1.9,0.46) node {$(p+1)_p$};
\fill [color=qqqqff] (-1.3,-0.68) circle (1.5pt);
\draw[color=qqqqff] (-1.4,-0.9) node {$(p+1)_2$};
\fill [color=qqqqff] (-0.66,-0.92) circle (1.5pt);
\draw[color=qqqqff] (-0.55,-1.2) node {$(p+1)_1$};
\end{scriptsize}
\end{tikzpicture}
\end{picture}
\label{fig1}
}
\qquad
\subfloat[The Graph $\mathcal{M}_2$.]{
\begin{picture}(200,200)(75,75)
\definecolor{qqqqff}{rgb}{0,0,0}
\begin{tikzpicture}[line cap=round,line join=round,>=triangle 45,x=1.0cm,y=1.0cm]
\clip(-5.66,-4.54) rectangle (15.68,4.94);
\draw [shift={(0.55,1.41)}] plot[domain=0.05:3.86,variable=\t]({1*1.61*cos(\t r)+0*1.61*sin(\t r)},{0*1.61*cos(\t r)+1*1.61*sin(\t r)});
\draw [shift={(0.56,1.38)},dash pattern=on 2pt off 2pt] plot[domain=-2.44:0.07,variable=\t]({1*1.6*cos(\t r)+0*1.6*sin(\t r)},{0*1.6*cos(\t r)+1*1.6*sin(\t r)});
\draw [shift={(1.86,1.62)},dash pattern=on 2pt off 2pt] plot[domain=-0.45:0.33,variable=\t]({1*1.39*cos(\t r)+0*1.39*sin(\t r)},{0*1.39*cos(\t r)+1*1.39*sin(\t r)});
\draw [shift={(-0.42,0.53)},dash pattern=on 2pt off 2pt] plot[domain=3.31:4.09,variable=\t]({1*1.48*cos(\t r)+0*1.48*sin(\t r)},{0*1.48*cos(\t r)+1*1.48*sin(\t r)});
\draw [shift={(-0.51,2.06)},dash pattern=on 2pt off 2pt] plot[domain=2.06:2.88,variable=\t]({1*1.45*cos(\t r)+0*1.45*sin(\t r)},{0*1.45*cos(\t r)+1*1.45*sin(\t r)});
\draw (-0.78,2.33)-- (-1.98,1.9);
\draw (-1.91,2.44)-- (-0.78,2.33);
\draw (-1.19,3.34)-- (-0.78,2.33);
\draw (1.2,4.)-- (0.88,3);
\draw (-0.66,-0.92)-- (-0.67,0.36);
\draw (-1.29,-0.68)-- (-0.67,0.36);
\draw (-1.87,0.28)-- (-0.67,0.36);
\draw (0.88,3)-- (-0.67,0.36);
\draw (3.12,1.02)-- (2.16,1.49);
\draw (3.17,2.08)-- (2.16,1.49);
\draw (2.78,2.56)-- (2.16,1.49);
\draw (-0.78,2.33)-- (2.16,1.49);
\draw (2.16,1.49)-- (-0.67,0.36);
\begin{scriptsize}
\fill [color=qqqqff] (-0.78,2.33) circle (1.5pt);
\draw[color=qqqqff] (-0.67,2.1) node {$1$};
\fill [color=qqqqff] (-1.19,3.34) circle (1.5pt);
\draw[color=qqqqff] (-1.17,3.6) node {$1_p$};
\fill [color=qqqqff] (-1.91,2.44) circle (1.5pt);
\draw[color=qqqqff] (-2.15,2.55) node {$1_2$};
\fill [color=qqqqff] (-1.98,1.9) circle (1.5pt);
\draw[color=qqqqff] (-2.2,2.0) node {$1_1$};
\fill [color=qqqqff] (0.88,3) circle (1.5pt);
\draw[color=qqqqff] (0.9,2.8) node {$2$};
\fill [color=qqqqff] (1.2,4.) circle (1.5pt);
\draw[color=qqqqff] (1.2,4.2) node {$2_1$};
\fill [color=qqqqff] (2.16,1.49) circle (1.5pt);
\draw[color=qqqqff] (2.05,1.3) node {$3$};
\fill [color=qqqqff] (3.12,1.02) circle (1.5pt);
\draw[color=qqqqff] (3.35,.9) node {$3_p$};
\fill [color=qqqqff] (3.17,2.08) circle (1.5pt);
\draw[color=qqqqff] (3.37,2.2) node {$3_2$};
\fill [color=qqqqff] (2.78,2.56) circle (1.5pt);
\draw[color=qqqqff] (2.89,2.74) node {$3_1$};
\fill [color=qqqqff] (-0.67,0.36) circle (1.5pt);
\draw[color=qqqqff] (-0.2,0.3) node {$p+1$};
\fill [color=qqqqff] (-1.87,0.28) circle (1.5pt);
\draw[color=qqqqff] (-1.9,0.46) node {$(p+1)_p$};
\fill [color=qqqqff] (-1.3,-0.68) circle (1.5pt);
\draw[color=qqqqff] (-1.4,-0.9) node {$(p+1)_2$};
\fill [color=qqqqff] (-0.66,-0.92) circle (1.5pt);
\draw[color=qqqqff] (-0.55,-1.2) node {$(p+1)_1$};
\end{scriptsize}
\end{tikzpicture}
\end{picture}
\label{fig2}
}
\caption{The Graphs of $\mathcal{M}_1$ and $\mathcal{M}_2$.}
\end{figure}

\begin{thm} \label{t2}
Suppose $p$ is a prime number and $G$ is a group with center $Z$ such that $\frac{G}{Z} \cong \mathbb{Z}_{p^2} \rtimes \mathbb{Z}_{p^2}$. The commuting conjugacy class graph of $G$ has one of the following types:
\begin{enumerate}[$(1)$]
\item \textit{$\frac{G}{Z}$ is abelian}. In this case, $\frac{G}{Z}\cong \mathbb{Z}_{p^2}\times \mathbb{Z}_{p^2}$ and the commuting conjugacy class graph of $G$ is isomorphic to $\mathcal{M}_1[K_{n(p-1)},\cdots,K_{n(p-1)},K_{m(p^2-p)},\cdots,K_{m(p^2-p)}]$, where there are $p+1$ copies of $K_{n(p-1)}$, $p^2+p$ copies of $K_{m(p^2-p)}$ and $\mathcal{M}_1$ is the graph depicted in Figure \ref{fig1}. Here, $m = \frac{|Z|}{p^2}$ and $n = \frac{|Z|}{p}$.

\item \textit{$\frac{G}{Z}$ is not abelian}. In this case, the commuting conjugacy class graph of $G$ is isomorphic to $\mathcal{M}_2[K_{n(p-1)},\cdots,K_{n(p-1)},  K_{n(p^2-p)}]$, where there are $p^2 + p + 1$ copies of $K_{n(p-1)}$, a copy of $K_{n(p^2-p)}$ and $\mathcal{M}_2$ is depicted in Figure \ref{fig2}. Here, $p$ is an odd prime and $n = \frac{|Z|}{p}$.
\end{enumerate}
\end{thm}

\label{fig}



\section{Preliminary Results}
The aim of this section is to prove some preliminary results which are crucial in next sections. The greatest common divisor of positive  integers $r$ and $s$ is denoted by $(r,s)$. We start this section by the following simple lemma:

\begin{lem}\label{lem3.2}
Let $p$ be a prime number and $G$ be a group with center $Z$ such that $x^{p^2},y^p\in Z$. Then, for each positive integer $m$ with this condition that $m \not\equiv 0 \ (mod \ p)$, we have $C_G(x^m)=C_G(x)$ and $C_G(y^m)=C_G(y)$.
\end{lem}

\begin{proof}
By our assumption, $(m,p^2)=1$ and there are $n$ and $k$ such that  $mn+kp^2=1$. Since $x^{p^2}\in Z$, $C_G(x^m)$ $\subseteq$ $C_G(x^{mn})$ $=$ $C_G(x^{1-kp^2})$ $=$ $C_G(xz)=C_G(x)$, as desired. The proof of second one is similar and so it is omitted.
\end{proof}

\begin{lem}\label{lem3.2.5}
Let $G$ be a finite group with center $Z$, $p \mid |G|$ and  $a,b\in G\setminus Z$ such that $o(aZ)=o(bZ)=p^2$ and $G$ $=$ $\{ a^ib^jz \mid z \in Z\}$. Then the following hold:

\begin{enumerate}[$(a)$]
\item  Suppose $1\leq i,j\leq p^2-1$ and at least one of $i, j$ is not divisible by $p$. Then $a^ib^j\neq b^ja^i$.

\item   $C_G(a^sb)=C_G(a^tb)$, $1\leq s,t\leq p^2-1$, if and only if $s=t$.
\end{enumerate}
\end{lem}

\begin{proof} Our argument consists in showing:
\begin{enumerate}[$(a)$]
\item \textit{Suppose $1\leq i,j\leq p^2-1$, $a^ib^j = b^ja^i$ and at least one of $i, j$ is not divisible by $p$}. Without loss of generality we can assume that $j \not\equiv 0 \ (mod \ p)$.  Then $(j,p^2)=1$ and so there are integers  $m$ and $n$ such that $mj+np^2 = 1$. Since  $b^{p^2}\in Z$, $a^ib = ba^i$ which means that $a^i\in Z$. This implies that $p^2=o(aZ)\mid i$, contradict by our condition that $1\leq i\leq p^2-1$.

\item If $s \ne t$ and $C_G(a^sb)=C_G(a^tb)$ then $a^sb\in C_G(a^tb)$ and so $a^{s-t}b  =  ba^{s-t}$. Again $p^2=o(aZ)\mid s-t$ which is impossible.
\end{enumerate}
Hence the result.
\end{proof}

The following result is a particular case of \cite[Theorem 8.1(c)]{maj}.

\begin{thm}\label{th2.0.0}
The non-abelian split extension $\mathbb{Z}_4\rtimes \mathbb{Z}_4$ is not  capable.
\end{thm}

Suppose $p$ is a prime number. Then it is easy to see that there are only two groups of order $p^4$ that can be written as a semidirect product of two cyclic groups of order $p^2$. One of these groups is the abelian group $L_1=\mathbb{Z}_{p^2} \times \mathbb{Z}_{p^2}$ and another one is a group presented as 
$L_2$ $=$ $\langle x,y\ \mid \ x^{p^2}=y^{p^2}=1, yx=x^{p+1}y\rangle.$, see \cite[pp. 145,146]{3} for details.
For simplicity of our argument, it is useful to write the presentations of $L_1$ and $L_2$ in the  form $\langle x,y\ \mid \ x^{p^2}=y^{p^2}=1, yx=x^{rp+1}y\rangle$, where $r = 0, 1$.

\begin{lem}\label{lem1}
Suppose $G = L_1$ or $L_2$. Then  for every $i,j$, $0\leq i,j\leq p^2-1$, which are not  simultaneously zero, we have:
\[o(x^iy^j)=
\begin{cases}
p &  p\mid (i,j)\\
p^2 & \text{Otherwise}
\end{cases}.\]
\end{lem}

\begin{proof}
The result   is trivial for $L_1$. So, it is enough to assume that $r=1$.  By definition, for every integers $i$ and $j$,  $0\leq i,j\leq p^2-1$, such that $i, j$ are not  simultaneously zero, we have
\begin{eqnarray*}
yx^2 & = & yxx \\
       & = & (x^{p+1}y)x=x^{p+1}(yx) \\
       & = & x^{p+1}(x^{p+1}y) \\
       & = & x^{2(p+1)}y.
\end{eqnarray*}
By the same way for each $i$, $yx^i=x^{i(p+1)}y$ and so
\begin{eqnarray*}
y^2x^i & = &  y(yx^i) \\
          & = & y(x^{i(p+1)}y)=(yx^{i(p+1)})y \\
          & = & (x^{i(p+1)^2}y)y=x^{i(p+1)^2}y^2.
\end{eqnarray*}
On the other hand, $y^jx^i=x^{i(p+1)^j}y^j$. Since $x^{p^2}=1$ and $(p+1)^j$ $=$ $p^j$ $+$ $jp^{j-1}$ $+$ $\cdots+\frac{j(j-1)}{2}p^2+jp+1$. So,
\begin{equation}\label{e1}
y^jx^i=x^{ijp+i}y^j.
\end{equation}
This shows that
\begin{eqnarray*}
(x^iy^j)^2 & = & (x^iy^j)(x^iy^j)=x^i(y^jx^i)y^j \\
               & = & x^i(x^{ijp+i}y^j)y^j\\
               & = & x^{ijp}x^{2i}y^{2j}.
\end{eqnarray*}
An inductive argument now proves that, for each $k$,
\begin{equation}\label{e2}
(x^iy^j)^k=x^{\frac{k(k-1)}{2}ijp}x^{ki}y^{kj}.
\end{equation}
Therefore, for every $i,j$ with  $0\leq i,j\leq p^2-1$ and $i,j$ are not simultaneously zero, we have:
\[o(x^iy^j)=
\begin{cases}
p &  p\mid (i,j)\\
p^2 & \text{otherwise}
\end{cases}.\]
This completes the proof.
\end{proof}

\section{Proof of Theorem 1.1}

In this section, the proof of our first main result will be presented. To do this, we need some information about a group $G$ with this property that  $|\frac{G}{Z(G)}| = p^4$, $p$ is prime, and $\frac{G}{Z(G)}$ can be generated by two elements of order $p^2$.

\begin{lem}\label{lem2.4}
Let $G$ be a finite group with center $Z$, $p$ be  prime  and  $a,b\in G\setminus Z$. Moreover, we assume that  $a^{p^2},b^{p^2}\in Z$  and $(a^ib^j)^k=a^{\frac{k(k-1)}{2} ijrp}a^{ki}b^{kj}z$ where $z\in Z$, $r=0, 1$ and $i,j,k$ are positive integers. The following statements hold:

\begin{enumerate}[$i)$]
\item If $1\leq i,j\leq p-1$, then  there are $n_j$ and $s$ such that $a^{ip}b^{jp}=(a^{sp}b^p)^{n_j}z$, where $n_j \not\equiv 0 \ (mid \ p)$, $1\leq s\leq p-1$ and $z\in Z$.

\item If $1\leq i\leq p-1$, $1\leq j\leq p^2-1$ and $j \not\equiv 0 \ (mod \ p)$, then  there are $n_j$ and $s$ such that $a^{ip}b^j=(a^{sp}b)^{n_j}z$, where $n_j \not\equiv 0 \ (mod \ p)$, $1\leq s\leq p-1$ and $z\in Z$.

\item If $1\leq j\leq p-1$, $1\leq i\leq p^2-1$ and $i \not\equiv 0 \ (mod \ p)$, then  there are $n_i$ and $s$ such that $a^ib^{jp}=(ab^{sp})^{n_i}z$, where $n_i \not\equiv 0 \ (mod \ p)$, $1\leq s\leq p-1$ and $z\in Z$.

\item If $1\leq i,j\leq p^2-1$, $i \not\equiv 0 \ (mod \ p)$ and $j \not\equiv 0 \ (mod \ p)$, then for  there are $n_j$ and $s$ such that $a^ib^j=(a^sb)^{n_j}z$, where  $n_j \not\equiv 0 \ (mod \ p)$, $s \not\equiv 0 \ (mod \ p)$, $1\leq s\leq p^2-1$ and $z\in Z$.

\item If $1\leq i,j\leq p^2-1$, $i \not\equiv 0 \ (mod \ p)$ and $j \not\equiv 0 \ (mod \ p)$, then there are $n_i$ and $s$ such that $a^ib^j=(ab^s)^{n_i}z$, where $n_i \not\equiv 0 \ (mod \ p)$, $s \not\equiv 0 \ (mod \ p)$, $1\leq s\leq p^2-1$ and $z\in Z$.
\end{enumerate}
\end{lem}

\begin{proof} By assumption, $(a^ib^j)^k=a^{\frac{k(k-1)}{2} ijrp}a^{ki}b^{kj}z$. In this equality, we assume that $i=k=p$ and $j=1$. Since $a^{p^2}\in Z$, $(a^pb)^p=b^pz'$ in which $z'\in Z$. Also $(a^pb)^p(a^pb)=(a^pb)(a^pb)^p$ then $(b^pz')(a^pb)=(a^pb)(b^pz')$ so $a^pb^p=b^pa^p$. Since $a^{p^2},a^{p^2}\in Z$,
\[(a^{ip}b^{jp})^p=a^{\frac{p(p-1)}{2} ijrp^3}a^{ip^2}b^{jp^2}z\in Z\]
then $(a^{ip}b^{jp})^p\in Z$. Also we can show that $(a^{ip}b^j)^{p^2}, (a^ib^j)^{p^2}\in Z$. Now, each part of the lemma will be proved separately.
\begin{enumerate}[$i)$]
\item In this case  $(j,p)=1$ and so  there are integers $m_j$ and $k$ such that $m_jj+kp=1$. By our assumption, $a^pb^p=b^pa^p$ and $b^{p^2}\in Z$ and hence,
\begin{eqnarray*}
(a^{ip}b^{jp})^{m_j} & = & a^{im_jp}b^{jm_jp} \\
                          & = & a^{im_jp}b^{(1-kp)p} \\
                          & = & a^{im_jp}b^pz_1.
\end{eqnarray*}
Since  $(a^{ip}b^{jp})^p\in Z$, $a^{ip}b^{jp}  =  (a^{im_jp}b^p)^{n_j}z_3$. It is clear that $p \nmid im_j$ and so  there are integers $s$ and $k$ such that $im_j=kp+s$. It is now easy to see that for a fixed positive integer $j$, $i_1$ $=$ $i_2$ if and only if $s_1=s_2$. Since $1\leq i\leq p-1$, $1\leq s\leq p-1$. So, by  above discussion and this fact that $a^{p^2}\in Z$,
\[a^{ip}b^{jp}=(a^{im_jp}b^p)^{n_j}z_3=(a^{kp^2+sp}b^p)^{n_j}z_3=(a^{sp}b^p)^{n_j}z,\]
as desired.

\item Since $p \nmid j$, $(j,p^2)=1$. So, there are $m_j$ and $k$ such that $m_jj+kp^2=1$ in which $(m_j,p)=1$. By our assumption $a^{p^2},b^{p^2}\in Z$ and $(a^ib^j)^k=a^{\frac{k(k-1)}{2} ijtp}a^{ki}b^{kj}$. Thus,
\begin{eqnarray*}
(a^{ip}b^j)^{m_j} & = & a^{\frac{m_j(m_j-1)}{2} ijtp^2}a^{im_jp}b^{jm_j}z_1 \\
                 & = & a^{im_jp}b^{1-kp^2}z_2 \\
                 & = & a^{im_jp}bz_3.
\end{eqnarray*}
Again since $(m_j,p)=1$, $(m_j,p^2)=1$ and so there are $n_j$ and $t$ such that $n_jm_j+~tp^2=1$ in which $p \nmid n_j$. We now apply above discussion and this fact that  $(a^{ip}b^j)^{p^2}\in Z$ to derive the equality
$a^{ip}b^j  =  (a^{im_jp}b)^{n_j}z_5$. A similar argument shows that  there exists an integer $s$ such that $1\leq s\leq p-1$ and $a^{ip}b^j=(a^{sp}b)^{n_j}z$ which completes our argument.

\item The proof is similar to $(ii)$ and so it is omitted.

\item Since $p \nmid j$, $(j,p^2)=1$. So there are $m_j$ and $k$ such that $m_jj+kp^2=1$ in which $(m_j,p)=1$. By our assumption $b^{p^2}\in Z$ and $(a^ib^j)^k=a^{\frac{k(k-1)}{2} ijtp}a^{ki}b^{kj}$. Therefore,
\begin{eqnarray*}
(a^ib^j)^{m_j} & = & a^{\frac{m_j(m_j-1)}{2} ijtp}a^{im_j}b^{jm_j}z_1 \\
                 & = & a^{\frac{m_j(m_j-1)}{2} ijtp+im_j}b^{1-kp^2}z_1 \\
                 & = & a^ubz_2,
\end{eqnarray*}
where
\begin{equation}\label{e111-5}
u=\frac{m_j(m_j-1)}{2}ijtp+im_j=i\biggl(\frac{m_j(m_j-1)}{2}jtp+m_j\biggr).
\end{equation}
Since $(m_j,p)=1$, $(m_j,p^2)=1$ and again there are $n_j$ and $t$ such that $m_jn_j+tp^2=1$ in which $p \nmid n_j$. By above discussion and this fact that $(a^ib^j)^{p^2}\in Z$, $a^ib^j  =  (a^ub)^{n_j}z_4$.
It is now clear that $p \nmid u$ and so $p^2 \nmid u$. Choose integers $s$ and $k$ such that $u=kp^2+s$ and $p \nmid s$. By Equation \ref{e111-5}, $j$ is fixed, one can see that $i_1=i_2$ if and only if $u_1=u_2$ if and only if $s_1=s_2$.  Since $1\leq i\leq p^2-1$,  $1\leq s\leq p^2-1$ and by above discussion and the fact that $a^{p^2}\in Z$, $a^ib^j$ $=$ $(a^ub)^{n_j}z_4$ $=$ $(a^{kp^2+s}b)^{n_j}z_4$ $=$ $(a^sb)^{n_j}z$ which completes the proof of this part.

\item The proof is similar to $iv$ and so it is omitted.
\end{enumerate}
Hence the result.
\end{proof}

Suppose $p$ is a prime number. By a result of Baer \cite{1.5}, the abelian group $\mathbb{Z}_{p^2} \times \mathbb{Z}_{p^2}$ is capable. If $p$ is an odd prime, then the non-abelian group $\mathbb{Z}_{p^2}\rtimes \mathbb{Z}_{p^2}$ is also capable \cite{8}. This means that there exists a group  $G$  such that $\frac{G}{Z(G)}\cong \mathbb{Z}_{p^2}\times \mathbb{Z}_{p^2}$ or $\mathbb{Z}_{p^2}\rtimes \mathbb{Z}_{p^2}$. In the next theorem, we obtain the number and structure of the centralizers of $G$. We also proved in  Theorem \ref{th2.0.0}  that the non-abelian group $\mathbb{Z}_{4}\rtimes \mathbb{Z}_{4}$ is not capable. In the next theorem, if $p=2$, then we will assume that  $\frac{G}{Z}\cong \mathbb{Z}_{4}\times \mathbb{Z}_{4}$.

\begin{thm}\label{thm3.5}
Suppose $p$ is a prime number and $G$ is a group with center $Z$ such that $\frac{G}{Z} \cong \mathbb{Z}_{p^2} \rtimes \mathbb{Z}_{p^2}$. Then the number of centralizers of $G$ is $(p+1)^2+1$.
\end{thm}

\begin{proof}
By our assumption, it is enough to assume that  $G$ is non-abelian and $Z\neq 1$. By our discussion before Lemma \ref{lem1}, we can see that
\[\frac{G}{Z}\cong Z_{p^2}\rtimes Z_{p^2}=\langle x,y\ \mid\ x^{p^2}=y^{p^2}=1\ ,\ yx=x^{rp+1}y\rangle.\]
Thus, there are $a, b\in G\setminus Z$ such that
\[\frac{G}{Z}=\langle aZ,bZ\ \mid\ (aZ)^{p^2}=(bZ)^{p^2}=1\ ,\ (bZ)(aZ)=(aZ)^{rp+1}(bZ)\rangle\]
and so
\[G=\{ a^ib^jz\ \mid\ 0\leq i,j\leq p^2-1\ ,\ ba=a^{rp+1}bz'\ ,\ z,z',a^{p^2},b^{p^2}\in Z\}.\]
We now apply Equation \ref{e2} to prove that
\begin{equation}\label{e33}
(a^ib^j)^k=a^{\frac{k(k-1)}{2}ijrp}a^{ki}b^{kj}z,
\end{equation}
where $z\in Z$. In Equation \ref{e33}, set $i=k=p$ and $j=1$. Since $a^{p^2}\in Z$, there exists $z\in Z$ such that $(a^pb)^p=b^pz$. Therefore,
\begin{eqnarray}
a^pb^p            & = & b^pa^p.  \label{e45}
\end{eqnarray}
So, for every $0\leq i,j\leq p-1$,
\begin{equation}\label{e5}
\{a^{sp}b^{tp}z \mid 0\leq s,t\leq p-1,\ z\in Z\}\subseteq C_G(a^{ip}b^{jp}z).
\end{equation}
Note that $|G|=p^4|Z|$ and for each $x\in G\setminus Z$, $Z\lneqq C_G(x)\lneqq G$. So,
\begin{equation}\label{e9}
|C_G(x)|=p|Z|\quad \text{or} \quad p^2|Z|\quad \text{or} \quad p^3|Z|.
\end{equation}
Suppose $i$ and $j$ are fixed and $0\leq i,j\leq p^2-1$. By Equation \ref{e33}, $(a^{i}b^{j})^{p^2}\in Z$. Thus,
\begin{equation}\label{e10}
\{(a^{i}b^{j})^kz \mid 0\leq k\leq p^2-1,\ z\in Z\}\subseteq C_G(a^{i}b^{j}).
\end{equation}

Next we partition the elements of $G$ into the following eight parts:
\begin{eqnarray*}
V_1 &=& \{ a^{i}z \mid p \mid i \ \& \ z\in Z\},\\
V_2 &=& \{ a^{i}z \mid p \nmid i \ \& \ z\in Z\},\\
V_3 &=& \{ b^{j}z \mid p \mid j \ \& \ z\in Z\},\\
V_4 &=& \{ b^{j}z \mid p \nmid j \ \& \ z\in Z\},\\
V_5 &=& \{ a^{i}b^{j}z \mid p \mid (i,j)\ \& \ z\in Z\},\\
V_6 &=& \{ a^{i}b^{j}z \mid p \mid i, \ p \nmid j \ \& \ z\in Z\},\\
V_7 &=& \{ a^{i}b^{j}z \mid p \nmid i, \ p \mid j \ \& \ z\in Z\},\\
V_8 &=& \{ a^{i}b^{j}z \mid p \nmid i, \ p \nmid j \ \& \ z\in Z\}.
\end{eqnarray*}

To compute the centralizer of non-central elements of $G$, the following cases will be considered:
\begin{enumerate}[a)]
\item \textit{$i\neq 0$ and $j=0$}. Then for every $1\leq i\leq p^2-1$,
\begin{equation}\label{e11}
\{a^{k}z\ \mid\ 0\leq k\leq p^2-1\ ,\ z\in Z\}\subseteq C(a^i).
\end{equation}

\begin{enumerate}[1)]
\item $p \mid i$. Suppose $i=sp$ such that $1\leq s\leq p-1$. By Equation \ref{e5},
\[A=\{a^kb^{tp}z \mid 0\leq k\leq p^2-1,\ 0\leq t\leq p-1,\ z\in Z\}\subseteq C_G(a^{i})=C_G(a^{sp}).\]
Since $|A|=p^3|Z|$,  by Equation \ref{e9} and Lemma \ref{lem3.2},
\[C_G(a^p)=C_G(a^{sp})=\{a^kb^{tp}z \mid 0\leq k\leq p^2-1,\ 0\leq t\leq p-1,\ z\in Z\}.\]

\item $p \nmid i$. By Lemma \ref{lem3.2.5}($a$), the element $a$  commutes only with elements of $V_1 \cup V_2$ and by Equations \ref{e11}, \ref{e9} and Lemma \ref{lem3.2},
\[C_G(a)=C_G(a^i)=\{a^kz \mid 0\leq k\leq p^2-1,\ z\in Z\}\leq C(a^p).\]
\end{enumerate}

\item \textit{$j\neq 0$ and $i=0$}. Apply a similar argument as the case of $(1)$. To do this, we consider two cases that $p \mid j$ or $p \nmid j$.
\begin{enumerate}[1)]
\setcounter{enumii}{2}
\item $p \mid j$. Suppose $j=tp$, $1\leq t\leq p-1$. Then,
\[C_G(b^p)=C_G(b^{tp})=\{a^{sp}b^kz \mid 0\leq k\leq p^2-1,\ 0\leq s\leq p-1,\ z\in Z\}.\]
\item $p \nmid j$. In this case, we can easily see that,
\[C_G(b)=C_G(b^j)=\{b^kz \mid 0\leq k\leq p^2-1,\ z\in Z\}\leq C(b^p).\]
\end{enumerate}

\item \textit{$i\neq 0$ and $j\neq 0$}. The following four subcases are considered into account.

\begin{enumerate}[1)]
\setcounter{enumii}{4}
\item \textit{$p \mid i$ and $p \mid j$}. Suppose $i=i_1p$ and $j=j_1p$ in which $1\leq i_1,j_1\leq p-1$.               By Lemma \ref{lem2.4}(I), there are $n_{j_1}$ and $s$ such that $a^{i_1p}b^{j_1p}=(a^{sp}b^p)^{n_{j_1}}z$, where $p \nmid n_{j_1}$, $1\leq s\leq p-1$ and $z\in Z$. By Equation \ref{e33}, $(a^{sp}b^p)^p\in Z$ and by Lemma (\ref{lem3.2}),
\[C_G(a^{i_1p}b^{j_1p})=C_G((a^{sp}b^p)^{n_{j_1}}z)=C_G((a^{sp}b^p)^{n_{j_1}})=C_G(a^{sp}b^p).\]
Since $1\leq s\leq p-1$,  the number of distinct proper centralizers constructed from the elements of $V_5$ is at most $p-1$. We now obtain the structures of these centralizers. In Equation \ref{e33}, we put  $i=s$, $j=1$ and $k=p$. Then $(a^sb)^p=a^{sp}b^pz$ in which $z\in Z$. Thus, $(a^{sp}b^p)(a^sb)$ $=$ $(a^sb)(a^{sp}b^p)$ and by Equation \ref{e5},
\[A=\{a^{ip}(a^sb)^kb^{jp}z \mid 0\leq i,j,k\leq p-1,\ z\in Z\}\subseteq C_G(a^{sp}b^p).\]
It can easily see that $|A|=p^3|Z|$ and by applying Equations \ref{e33} and \ref{e9},
\begin{equation}\label{e12}
C_G(a^{sp}b^p)=\{a^{ip}a^{sk}b^kb^{jp}z \mid 0\leq i,j,k\leq p-1,\ z\in Z\}.
\end{equation}
Suppose $1\leq t\leq p-1$. By Equation \ref{e12},
\begin{equation}\label{e12-5}
a^tb\in C_G(a^{sp}b^p) \ \text{if \ and \ only \ if} \ t=s.
\end{equation}
Therefore,  $C_G(a^{tp}b^p)=C_G(a^{sp}b^p)$, $1\leq s,t\leq p-1$, if and only if $t=s$. This shows that the number of distinct proper centralizers is exactly $p-1$.

\item \textit{$p \mid i$ and $p \nmid j$}. Suppose $i=i_1p$, where $1\leq i_1\leq p-1$. By Lemma \ref{lem2.4}(II), there are $n_j$ and $s$ such that $a^{i_1p}b^j=(a^{sp}b)^{n_j}z$ in which $p \nmid n_j$, $1\leq s\leq p-1$ and $z\in Z$. Now by Equation \ref{e33}, $(a^{sp}b)^{p^2}\in Z$. Hence by Lemma \ref{lem3.2}, $C_G(a^{i_1p}b^j)$ $=$ $C_G((a^{sp}b)^{n_j}z)$ $=$ $C_G((a^{sp}b)^{n_j})$ $=$ $C_G(a^{sp}b)$. Since $1\leq s\leq p-1$, by Lemma \ref{lem3.2.5}(b), the number of distinct proper centralizers is exactly $p-1$. To obtain the structures of these centralizers, it is enough to apply Equation \ref{e10}. This implies that \[A=\{(a^{sp}b)^kz\ \mid\ 0\leq k\leq p^2-1\ ,\ z\in Z\}\subseteq C(a^{sp}b).\]
On the other hand, by Equation \ref{e33}, $(a^{sp}b)^k=a^{spk}b^kz'$ in which $z'\in Z$. We assume that $k=tp+k_1$ in which $0\leq t,k_1\leq p-1$. Thus, $(a^{sp}b)^k$ $=$ $a^{spk}b^kz'$ $=$ $a^{spk_1}b^{k_1}b^{tp}z''$, where $z''\in Z$. So,
\begin{equation*}
A=\{a^{spk_1}b^{k_1}b^{tp}z\ \mid\ 0\leq k_1,t\leq p-1\ ,\ z\in Z\}\subseteq C(a^{sp}b).
\end{equation*}
By Lemma \ref{lem3.2.5}(a), $a^{sp}b$ cannot commute with an element of $V_1 \cup V_2 \cup V_4 \cup V_5 \cup V_7$. Furthermore, it is easy to see that $V_3 \subseteq A$. We claim that there are no elements $x \in V_6$ and $y \in V_8$ such that $xy=yx$. To see this, we assume that $x = a^{i_1p}b^{j_1} \in V_6$ and $y = a^{i_2}b^{j_2} \in V_8$ are arbitrary. By Lemma \ref{lem2.4}($iii$,$iv$), $a^{i_1p}b^{j_1}=(a^{sp}b)^{n_{j_1}}z_1$ and $a^{i_2}b^{j_2}=(a^tb)^{n_{j_2}}z_2$, where $p \nmid n_{j_1}$, $p \nmid n_{j_2}$, $p \nmid s$ and $p \nmid t$. By Equation \ref{e33}, $(a^{sp}b)^{p^2},(a^tb)^{p^2}\in Z$ and by
Lemma \ref{lem3.2}, $C_G(a^{i_1p}b^{j_1})=C_G(a^{sp}b)$ and $C_G(a^{i_2}b^{j_2})=C_G(a^tb)$. Thus, $(a^{i_1p}b^{j_1})(a^{i_2}b^{j_2})$ $=$ $(a^{i_2}b^{j_2})(a^{i_1p}b^{j_1})$ if and only if $(a^{sp}b)(a^tb)=(a^tb)(a^{sp}b)$ if and only if $a^{sp}ba^t=a^tba^{sp}$ if and only if $a^{sp-t}b=ba^{sp-t}$. We now apply Lemma \ref{lem3.2.5}(a) to deduce that $p\mid t$ which is impossible. Next we prove that $C(a^{sp}b)=A$. Choose $x\in C(a^{sp}b)\setminus A$. By above discussion, $x \in V_6$ and so we can write $x=a^{ip}b^{j}$, where $1\leq i\leq p-1$. By Lemma \ref{lem2.4}($ii$), there are $n_{j}$ and $m$ such that $a^{ip}b^{j}=(a^{mp}b)^{n_{j}}z$ with this condition that $p \nmid n_{j}$, $1\leq m\leq p-1$, and $z\in Z$. On the other hand, By Equation \ref{e33}, $(a^{mp}b)^{p^2}\in Z$. So, by Lemma \ref{lem3.2}, $C(a^{ip}b^j)=C(a^{mp}b)$. Since $x\in C(a^{sp}b)$, $a^{sp}b\in C_G(x)=C(a^{mp}b)$. Hence, $a^{(s-m)p}b$ $=$ $ba^{(s-m)p}$. Again by Lemma \ref{lem3.2.5}(a), $p\mid s-m$ and since $1\leq s,m\leq p-1$, $s=m$. This proves that $C(a^{sp}b)=A$. Therefore,
\[C(a^{sp}b)=\{a^{spk}b^{k}b^{tp}z\ \mid\ 0\leq k,t\leq p-1\ ,\ z\in Z\}.\]

\item \textit{$p \mid j$ and $p \nmid i$}. A similar argument as Case (6) shows that the number of distinct proper centralizers of $G$ is equal to $p-1$ and
\[C_G(ab^{tp})=\{a^{sp}a^kb^{tpk}z\mid 1\leq k,s\leq p-1,\ z\in Z\}.\]

\item \textit{$p \nmid i$ and $p \nmid j$}. By Lemma \ref{lem2.4}($iv$), there are $n_j$ and $s$ such that $a^ib^j=(a^sb)^{n_j}z$ in which $p \nmid n_j$, $p \nmid s$, $1\leq s\leq p^2-1$ and $z\in Z$. By Equation  \ref{e33}, $(a^sb)^{p^2}\in Z$. Thus, by Lemma \ref{lem3.2}, $C_G(a^ib^j)$ $=$ $C_G((a^sb)^{n_j}z)$ $=$ $C_G((a^sb)^{n_j})$ $=$ $C_G(a^sb)$. Since $p \nmid s$, $1\leq s\leq p^2-1$ and Lemma \ref{lem3.2.5}(b), the number of distinct proper centralizers contained in $V_8$ is equal to $p^2-p$. We now obtain the structure of centralizers. By Equations \ref{e10} and (\ref{e33}),
\begin{eqnarray}
A & = & \{(a^sb)^kz\ \mid\ 0\leq k\leq p^2-1\ ,\ z\in Z\}\subseteq C(a^sb) \notag \\
                         & = & \{a^{\frac{k(k-1)}{2}sp+sk}b^kz\ \mid\ 0\leq k\leq p^2-1\ ,\ z\in Z\}. \label{e20}
\end{eqnarray}
By Lemma \ref{lem3.2.5}(a) and definition of $V_7$, one can see that $a^sb$ cannot be commuted with an elements in $V_1 \cup V_2 \cup V_3 \cup V_4 \cup V_6 \cup V_7$. We claim that $C(a^{sp}b)=A$. To prove, we assume that $x\in C_G(a^sb)\setminus A$. By above discussion, $x \in V_5 \cup V_8$. If $x\in V_8$, then $x=a^{i'}b^{j'}$ and by Lemma \ref{lem2.4}($iv$), there are $n_{j'}$ and $t$ such that $a^{i'}b^{j'}=(a^tb)^{n_{j'}}z$ in which $p \nmid n_{j'}$, $p \nmid t$ and $z\in Z$. Hence by Equation \ref{e33}, $(a^tb)^{p^2}\in Z$ and by Lemma \ref{lem3.2}, we have $C_G(a^{i'}b^{j'})=C_G(a^tb)$. Since $x\in C_G(a^sb)$, $a^sb\in C_G(x)=C_G(a^tb)$ and so   $a^{s-t}b$  $=$  $ba^{s-t}$. Apply Lemma \ref{lem3.2.5}(a) to deduce that $p^2\mid s-t$. Since $1\leq s,t\leq p^2-1$, $s=t$ and $x$ $=$ $a^{i'}b^{j'}$ $=$ $(a^sb)^{n_{j'}}z\in A$. If $x\in V_5$, then $x=a^{i'p}b^{j'p}$ and by Lemma \ref{lem2.4}($i$), there are $n_{j'}$ and $t$ such that $a^{i'p}b^{j'p}=(a^{tp}b^p)^{n_{j'}}z$ in which $p \nmid n_{j'}$, $p \nmid t$ and $z\in Z$. Next by Equation \ref{e33}, $(a^{tp}b^p)^p\in Z$ and Lemma \ref{lem3.2} implies that $C_G(a^{i'p}b^{j'p})=C_G(a^{tp}b^p)$. Since $x\in C_G(a^sb)$, $a^sb\in C_G(x)=C_G(a^{tp}b^p)$ and  by Equation \ref{e12-5}, one can see that $t=s$. Hence Equation \ref{e33} shows that $x$ $=$ $a^{i'p}b^{j'p}$ $=$ $(a^{sp}b^p)^{n_{j'}}z$ $=$ $(a^sb)^{n_{j'}p}z\in A$. Therefore, in both cases $x\in A$ which means that
\[C_G(a^sb)=A=\{(a^sb)^kz\ \mid\ 0\leq k\leq p^2-1\ ,\ z\in Z\}.\]
In addition, by Equations \ref{e12} and \ref{e20}, we can  see that $C_G(a^sb)\leq C_G(a^{sp}b^p)$.
\end{enumerate}
\end{enumerate}

Therefore, the structure of non-trivial proper centralizers of $G$ is as follows:
\begin{eqnarray*}
C_G(a^p) & = & \{a^kb^{jp}z \mid 0\leq j\leq p-1,\ 0\leq k\leq p^2-1,\ z\in Z\},\\
C_G(a) & = & \{a^kz\mid 0\leq k\leq p^2-1\}\leq C_G(a^p),\\
C_G(b^p) & = & \{a^{ip}b^kz \mid 0\leq i\leq p-1,\ 0\leq k\leq p^2-1,\ z\in Z\},\\
C_G(b) & = & \{b^kz\mid 0\leq k\leq p^2-1\}\leq C_G(b^p),\\
C_G(a^{sp}b^p) & = & \{a^{ip}a^{sk}b^{k}b^{jp}z \mid 0\leq i,j,k\leq p-1,\ z\in Z\},  1\leq s\leq p-1\\
C_G(ab^{sp}) & = & \{a^{ip}a^kb^{spk}z \mid 0\leq i,k\leq p-1,\ z\in Z\}\leq C_G(a^{p}),  1\leq s\leq p-1\\
C_G(a^{sp}b) & = & \{a^{spk}b^kb^{jp}z \mid 0\leq j,k\leq p-1,\ z\in Z\}\leq C_G(b^{p}),  1\leq s\leq p-1\\
C_G(a^sb) & = & \{(a^sb)^ka^{isp}b^{ip}z\mid 0\leq i,k\leq p-1,\ z\in Z\}\leq C_G(a^{sp}b^p), 1\leq s\leq p^2-1,  p \nmid s
\end{eqnarray*}

So, the number of distinct proper centralizers of $G$ is  $(p^2-p)+3(p-1)+4=p^2+2p+1=(p+1)^2$ and therefore  $|Cent(G)|=(p+1)^2+1$. This completes the proof.
\end{proof}
We end this section by the following conjecture:

\begin{conj}
Suppose $p$ is a prime number, $n$ is a positive integer and $G$ is a group with center $Z$ such that $\frac{G}{Z}\cong \mathbb{Z}_{p^n}\rtimes \mathbb{Z}_{p^n}$. Then $|Cent(G)| = (p+1)^n+1$.
\end{conj}


\section{Proof of Theorem 1.2}


In this section, we apply the results of Section 3 to obtain the structure of the commuting conjugacy class graph of $G$ when $\frac{G}{Z(G)}\cong \mathbb{Z}_{p^2}\rtimes \mathbb{Z}_{p^2}$.


\begin{thm}
Suppose $p$ is a prime and $G$ is a group with center $Z$ such that $\frac{G}{Z}\cong \mathbb{Z}_{p^2}\times \mathbb{Z}_{p^2}$. Then the commuting conjugacy classes graph of $G$ has the graph structure
\[\mathcal{M}_1[\overbrace{K_{n(p-1)},\cdots,K_{n(p-1)}}^{p+1},\overbrace{K_{m(p^2-p)},\cdots,K_{m(p^2-p)}}^{p^2+p}],\]
where $\mathcal{M}_1$ is the graph depicted in Figure \ref{fig3}. Here, $m = \frac{|Z|}{p^2}$ and $n = \frac{|Z|}{p}$.

\vspace*{-2cm}
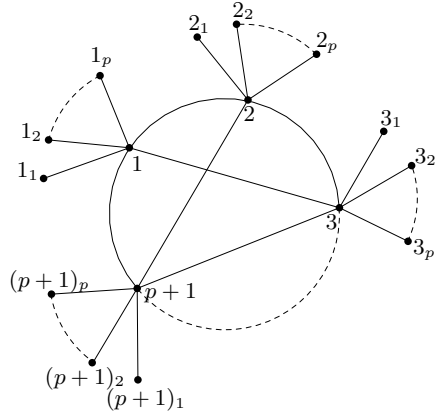
\begin{figure}[h]
\begin{picture}(370,300)
\scalebox{.95}{
\definecolor{qqqqff}{rgb}{0,0,0}
\begin{tikzpicture}[line cap=round,line join=round,>=triangle 45,x=1.0cm,y=1.0cm]
\clip(-5.66,-4.54) rectangle (15.68,4.94);
\draw [shift={(0.55,1.41)}] plot[domain=0.05:3.86,variable=\t]({1*1.61*cos(\t r)+0*1.61*sin(\t r)},{0*1.61*cos(\t r)+1*1.61*sin(\t r)});
\draw [shift={(0.56,1.38)},dash pattern=on 2pt off 2pt]  plot[domain=-2.44:0.07,variable=\t]({1*1.6*cos(\t r)+0*1.6*sin(\t r)},{0*1.6*cos(\t r)+1*1.6*sin(\t r)});
\draw [shift={(0.81,2.6)},dash pattern=on 2pt off 2pt]  plot[domain=0.79:1.63,variable=\t]({1*1.46*cos(\t r)+0*1.46*sin(\t r)},{0*1.46*cos(\t r)+1*1.46*sin(\t r)});
\draw [shift={(1.86,1.62)},dash pattern=on 2pt off 2pt]  plot[domain=-0.45:0.33,variable=\t]({1*1.39*cos(\t r)+0*1.39*sin(\t r)},{0*1.39*cos(\t r)+1*1.39*sin(\t r)});
\draw [shift={(-0.42,0.53)},dash pattern=on 2pt off 2pt]  plot[domain=3.31:4.09,variable=\t]({1*1.48*cos(\t r)+0*1.48*sin(\t r)},{0*1.48*cos(\t r)+1*1.48*sin(\t r)});
\draw [shift={(-0.51,2.06)},dash pattern=on 2pt off 2pt]  plot[domain=2.06:2.88,variable=\t]({1*1.45*cos(\t r)+0*1.45*sin(\t r)},{0*1.45*cos(\t r)+1*1.45*sin(\t r)});
\draw (-0.78,2.33)-- (-1.98,1.9);
\draw (-1.91,2.44)-- (-0.78,2.33);
\draw (-1.19,3.34)-- (-0.78,2.33);
\draw (0.17,3.88)-- (0.88,3);
\draw (0.72,4.06)-- (0.88,3);
\draw (1.84,3.64)-- (0.88,3);
\draw (-0.66,-0.92)-- (-0.67,0.36);
\draw (-1.29,-0.68)-- (-0.67,0.36);
\draw (-1.87,0.28)-- (-0.67,0.36);
\draw (0.88,3)-- (-0.67,0.36);
\draw (3.12,1.02)-- (2.16,1.49);
\draw (3.17,2.08)-- (2.16,1.49);
\draw (2.78,2.56)-- (2.16,1.49);
\draw (-0.78,2.33)-- (2.16,1.49);
\draw (2.16,1.49)-- (-0.67,0.36);
\begin{scriptsize}
\fill [color=qqqqff] (-0.78,2.33) circle (1.5pt);
\draw[color=qqqqff] (-0.67,2.1) node {$1$};
\fill [color=qqqqff] (-1.19,3.34) circle (1.5pt);
\draw[color=qqqqff] (-1.17,3.6) node {$1_p$};
\fill [color=qqqqff] (-1.91,2.44) circle (1.5pt);
\draw[color=qqqqff] (-2.15,2.55) node {$1_2$};
\fill [color=qqqqff] (-1.98,1.9) circle (1.5pt);
\draw[color=qqqqff] (-2.2,2.0) node {$1_1$};
\fill [color=qqqqff] (0.88,3) circle (1.5pt);
\draw[color=qqqqff] (0.9,2.8) node {$2$};
\fill [color=qqqqff] (1.84,3.64) circle (1.5pt);
\draw[color=qqqqff] (1.99,3.8) node {$2_p$};
\fill [color=qqqqff] (0.72,4.06) circle (1.5pt);
\draw[color=qqqqff] (0.82,4.25) node {$2_2$};
\fill [color=qqqqff] (0.17,3.88) circle (1.5pt);
\draw[color=qqqqff] (0.2,4.06) node {$2_1$};
\fill [color=qqqqff] (2.16,1.49) circle (1.5pt);
\draw[color=qqqqff] (2.05,1.3) node {$3$};
\fill [color=qqqqff] (3.12,1.02) circle (1.5pt);
\draw[color=qqqqff] (3.35,.9) node {$3_p$};
\fill [color=qqqqff] (3.17,2.08) circle (1.5pt);
\draw[color=qqqqff] (3.37,2.2) node {$3_2$};
\fill [color=qqqqff] (2.78,2.56) circle (1.5pt);
\draw[color=qqqqff] (2.89,2.74) node {$3_1$};
\fill [color=qqqqff] (-0.67,0.36) circle (1.5pt);
\draw[color=qqqqff] (-0.2,0.3) node {$p+1$};
\fill [color=qqqqff] (-1.87,0.28) circle (1.5pt);
\draw[color=qqqqff] (-1.9,0.46) node {$(p+1)_p$};
\fill [color=qqqqff] (-1.3,-0.68) circle (1.5pt);
\draw[color=qqqqff] (-1.4,-0.9) node {$(p+1)_2$};
\fill [color=qqqqff] (-0.66,-0.92) circle (1.5pt);
\draw[color=qqqqff] (-0.55,-1.2) node {$(p+1)_1$};
\end{scriptsize}
\end{tikzpicture}
}
\end{picture}
\vspace*{-3cm}
\caption{The Graph $\mathcal{M}_1$.}\label{fig3}
\end{figure}
\end{thm}

\begin{proof}
Since $\frac{G}{Z}$ is abelian, $xgZ=xZgZ=gZxZ=gxZ$, where $x,g\in G\setminus Z$ are arbitrary. Hence
there exists $z\in Z$ such $xg=gxz$ and so for every $x\in G\setminus Z$,
\begin{equation}\label{c.1}
x^G=\{g^{-1}xg \mid g\in G\}=\{g^{-1}gxz \mid z\in Z\}=\{xz \mid z\in Z\}.
\end{equation}

We now apply Theorem \ref{thm3.5} to prove that the group $G$ has the following structure:
\[G=\{a^ib^jz \mid 0\leq i,j\leq p^2-1,\ ba=abz_t,\ z,z_t\in Z\}.\]
Hence for each $x\in G\setminus Z$, $x=a^ib^jz$, where $0\leq i,j\leq p^2-1$ and $i, j$ are not simultaneously zero.
To obtain the non-central conjugacy classes of $G$, the following have to be investigated:
\begin{enumerate}[a)]
\item \textit{$i\neq 0$ and $j=0$.} We will consider two cases that $p \mid i$ and $p \nmid i$.
\begin{enumerate}[1)]
\item $p\mid i$. In this case, $i=sp$, where $1\leq s\leq p-1$. By Theorem \ref{thm3.5}, $|C_G(a^{sp})|=p^3|Z|$ and so $|(a^{sp})^G|=p$. By Equation \ref{c.1}, $(a^{sp})^G$ $=$ $\{a^{sp}z_1,a^{sp}z_2,\cdots,a^{sp}z_p\}$ $=$ $a^{sp}H$, where  $H=\{z_1,z_2,\cdots,z_p\}\subseteq Z$. Suppose $|Z|>p$ and choose $z_{r_1}\in Z\setminus H$. It is easy to see that $a^{sp}z_{r_1}\not\in (a^{sp})^G$ and so $(a^{sp}z_{r_1})^G\neq (a^{sp})^G$. Since $C_G(a^{sp}z_{r_1})=C_G(a^{sp})$, $|(a^{sp}z_{r_1})^G|=|(a^{sp})^G|=p$. Also,
$(a^{sp}z_{r_1})^G=(a^{sp})^Gz_{r_1}=a^{sp}Hz_{r_1}$. Since $(a^{sp})^G\cap(a^{sp}z_{r_1})^G=\emptyset$,
$H\cap Hz_{r_1}=\emptyset$ and $H\cup Hz_{r_1}\subseteq Z$. We now choose an element $z_{r_2}\in Z\setminus(H\cup Hz_{r_1})$ and by above method we will see that $p\mid |Z|$. Suppose $n=\frac{|Z|}{p}$. Thus, there are $n$ distinct conjugacy classes of the form $(a^{sp}z_{r_1})^G,(a^{sp}z_{r_2})^G, \ldots, (a^{sp}z_{r_n})^G$. Since $1\leq s\leq p-1$,  there are $n(p-1)$ distinct conjugacy classes with $n=\frac{|Z|}{p}$ and each of which has $p$ elements.

\item $p\nmid i$. Then by Theorem \ref{thm3.5}, $|C_G(a^i)|=p^2|Z|$ and so $|(a^i)^G|=p^2$. Also, by Equation \ref{c.1}, $(a^i)^G=\{a^iz_1,a^iz_2,\cdots,a^iz_{p^2}\}=a^iH$, where
$H=\{z_1,z_2,\cdots,z_{p^2}\}\subseteq Z$. Suppose that $|Z|>p^2$ and choose $z_{r_1}\in Z\setminus H$. It is easy to see that $a^iz_{r_1}\not\in (a^i)^G$ and hence $(a^iz_{r_1})^G\neq (a^i)^G$. Since $C_G(a^iz_{r_1})=C_G(a^i)$,  $|(a^iz_{r_1})^G|=|(a^i)^G|=p^2$. Also, $(a^iz_{r_1})^G=(a^i)^Gz_{r_1}=a^iHz_{r_1}$. Note that
$(a^i)^G\cap(a^iz_{r_1})^G=\emptyset$. Thus, $H\cap Hz_{r_1}=\emptyset$ and $H\cup Hz_{r_1}\subseteq Z$. Choose the element $z_{r_2}\in Z\setminus(H\cup Hz_{r_1})$. By repeated applications of above method we can see that $p^2\mid |Z|$. Define $m=\frac{|Z|}{p^2}$. Then for constant $i$, there are $m$ distinct conjugacy classes in the form of $(a^iz_{r_1})^G,(a^iz_{r_2})^G, \ldots, (a^iz_{r_m})^G$. Since $1\leq i\leq p^2-1$ and $p\nmid i$, there are $m(p^2-p)$ distinct conjugacy classes, where $m=\frac{|Z|}{p^2}$ and each of which has $p^2$ elements.
\end{enumerate}

\item \textit{$i=0$ and $j\neq0$.}
\begin{enumerate}[1)]
\setcounter{enumii}{2}
\item If $p\mid j$, then $j=tp$, where $1\leq t\leq p-1$. By Theorem \ref{thm3.5}, $|C_G(b^{tp})|=p^3|Z|$ and so $|(b^{tp})^G|=p$. Then  by Equation \ref{c.1} and a similar argument as (1),  for each $t$ there are $n$ distinct conjugacy classes of the form $(b^{tp}z_{r_1})^G,(b^{tp}z_{r_2})^G, \ldots, (b^{tp}z_{r_n})^G$, where $n=\frac{|Z|}{p}$.  Since $1\leq t\leq p-1$,  there are $n(p-1)$ distinct conjugacy classes with $n=\frac{|Z|}{p}$ and each class has $p$ elements.

\item If $p\nmid j$, then by Theorem \ref{thm3.5}, $|C_G(b^j)|=p^2|Z|$ and so
$|(b^j)^G|=p^2$. Apply again Equation \ref{c.1} and a similar argument as (2) to result that for each constant $j$, there are $m$ distinct conjugacy classes of the form $(b^jz_{r_1})^G,(b^jz_{r_2})^G, \ldots, (b^jz_{r_m})^G$, where $m=\frac{|Z|}{p^2}$. Since $1\leq j\leq p^2-1$ and
$p\nmid j$,  there are $m(p^2-p)$ distinct conjugacy classes with $m=\frac{|Z|}{p^2}$ and
each class has $p^2$ elements.
\end{enumerate}

\item $i\neq 0$ and $j\neq 0$. This case can be separated into the following four cases:
\begin{enumerate}[1)]
\setcounter{enumii}{4}

\item \textit{$p\mid i$ and $p\mid j$}. In this case, $i=sp$ and $j=tp$ in which $1\leq s,t\leq p-1$.
By Theorem \ref{thm3.5}, $|C_G(a^{sp}b^{tp})|=p^3|Z|$ and so $|(a^{sp}b^{tp})^G|=p$.
We now apply Equation \ref{c.1} and a similar argument as (1) for constants $s$ and $t$ to deduce that
there are $n$ distinct conjugacy classes of the form
$(a^{sp}b^{tp}z_{r_1})^G,(a^{sp}b^{tp}z_{r_2})^G, \ldots, (a^{sp}b^{tp}z_{r_n})^G$, where $n=\frac{|Z|}{p}$.
Since $1\leq s,t\leq p-1$,  there are $n(p-1)^2$ distinct conjugacy classes with $n=\frac{|Z|}{p}$
and each class has $p$ members.

\item \textit{$p\nmid i$ and $p\mid j$}. Clearly $j=tp$, $1\leq t\leq p-1$.
By Theorem \ref{thm3.5}, $|C_G(a^ib^{tp})|=p^2|Z|$ and hence $|(a^ib^{tp})^G|=p^2$. On the other hand, Equation \ref{c.1} and a similar argument as (2) for constants $i$ and $t$, show that there are $m$ distinct conjugacy classes of the form $(a^ib^{tp}z_{r_1})^G,(a^ib^{tp}z_{r_2})^G, \ldots, (a^ib^{tp}z_{r_m})^G$, where $m=\frac{|Z|}{p^2}$. Since $1\leq t\leq p-1$, $1\leq i\leq p^2-1$ and $p\nmid i$,  there are $mp(p-1)^2$ distinct conjugacy classes with $m=\frac{|Z|}{p^2}$ and each class has $p^2$ members.

\item \textit{$p\mid i$ and $p\nmid j$}. In this case $i=sp$, for some integer $s$ such that $1\leq s\leq p-1$.
By Theorem \ref{thm3.5}, $|C_G(a^{sp}b^j)|=p^2|Z|$ and so $|(a^{sp}b^j)^G|=p^2$. Then by Equation \ref{c.1} and a similar argument as (2) for constants $s$ and $j$, one can prove that there are $m$ distinct conjugacy classes of the  form $(a^{sp}b^jz_{r_1})^G,(a^{sp}b^jz_{r_2})^G, \ldots, (a^{sp}b^jz_{r_m})^G$, where $m=\frac{|Z|}{p^2}$. Since $1\leq s\leq p-1$, $1\leq j\leq p^2-1$ and $p\nmid j$,  there are $mp(p-1)^2$ distinct conjugacy classes with $m=\frac{|Z|}{p^2}$ and each class has $p^2$ members.

\item \textit{$p\nmid i$ and $p\nmid j$}. By Theorem \ref{thm3.5}, $|C_G(a^ib^j)|=p^2|Z|$ and
so $|(a^ib^j)^G|=p^2$. On the other hand, by Equation \ref{c.1} and applying a  similar argument as (2) for constants $i$ and $j$, one can see that there are $m$ distinct conjugacy classes of the form $(a^ib^jz_{r_1})^G,(a^ib^jz_{r_2})^G, \ldots, (a^ib^jz_{r_m})^G$, where $m=\frac{|Z|}{p^2}$. Since  $1\leq i,j\leq p^2-1$ and $p\nmid i,j$, there are $mp^2(p-1)^2$ distinct conjugacy classes with $m=\frac{|Z|}{p^2}$ and each class has length $p^2$.
\end{enumerate}
\end{enumerate}
The above discussion are summarized in Table \ref{j3.1}.
\begin{table}[h]
\begin{center}
\begin{tabular}{ |c|ll|c|l| }
\hline
 Type  & \multicolumn{2}{|c|}{ The Representatives of Conjugacy Classes }  & $|x^G|$ & $\#$ Conjugacy Classes \\ \hline
$1$ & $a^{sp}$ & $1\leq s\leq p-1$ & $p$ & $n(p-1)$ \\ \hline
$2$ & $a^i$      & $1\leq i\leq p^2-1,i\overset{p}{\not\equiv} 0$ & $p^2$ & $mp(p-1)$ \\ \hline
$3$ & $b^{tp}$ & $1\leq t\leq p-1$ & $p$ & $n(p-1)$ \\ \hline
$4$ & $b^j$      & $1\leq j\leq p^2-1,j\overset{p}{\not\equiv} 0$ & $p^2$  & $mp(p-1)$ \\ \hline
$5$ & $a^{sp}b^{tp}$ & $1\leq s,t\leq p-1$ & $p$  & $n(p-1)^2$ \\ \hline
$6$ & $a^ib^{tp}$ & $1\leq t\leq p-1,1\leq i\leq p^2-1,i\overset{p}{\not\equiv} 0$ & $p^2$  & $mp(p-1)^2$ \\ \hline
$7$ & $a^{sp}b^j$ & $1\leq s\leq p-1,1\leq j\leq p^2-1,j\overset{p}{\not\equiv} 0$ & $p^2$  & $mp(p-1)^2$ \\ \hline
$8$ & $a^ib^j$ & $1\leq i,j\leq p^2-1,i\overset{p}{\not\equiv} 0,j\overset{p}{\not\equiv} 0$ & $p^2$ & $mp^2(p-1)^2$ \\ \hline
\end{tabular}
\end{center}
\caption{Conjugacy classes of $G$ with $n=\frac{|Z|}{p}$ and $m=\frac{|Z|}{p^2}$.}\label{j3.1}
\end{table}

We now obtain the commuting conjugacy class graph of $G$. To do this, we consider the following cases:
\begin{enumerate}
\item Suppose $a^{s_1p}$ and $a^{s_2p}$ are representatives of two conjugacy classes  of Type 1 in Table 1.
It is clear that $a^{s_1p}a^{s_2p}=a^{s_2p}a^{s_1p}$ and so all such classes   are commuting together. Hence, the commuting conjugacy class graph has a subgraph isomorphic to the complete graph $K_{n(p-1)}$. We now assume that $a^{sp}$ is a representative of   a conjugacy class of Type 1  and $a^ub^vz$ is an arbitrary element of $G$ such that  $(a^{sp})(a^ub^vz) = (a^ub^vz)(a^{sp})$. Thus $a^{sp}b^v = b^va^{sp}$ and by Lemma \ref{lem3.2.5}(a), $v \equiv 0 \ (mod \ p)$ or $v=tp$ in which $0\leq t\leq p-1$. Hence $a^{sp}$ commutes with all elements in the form $a^ub^{tp}z$. Now the following cases can be occurred:
\begin{enumerate}
\item If $t=0$, then the conjugacy classes of Types 1 and 2 in Table 1 are commuting to each other.

\item If $t\neq 0$ and $u=0$, then the conjugacy classes of Types 1 and 3 in Table 1 are commuting to each other.

\item If $u,t\neq 0$ and $p\mid u$, then the conjugacy classes of Types 1 and 5 in Table 1 can be commuted to each other.

\item If $t\neq 0$ and $p\nmid u$, then the conjugacy classes of Types 1 and 6 in Table 1 are commuting together.
\end{enumerate}

\item It is clear that two conjugacy classes of Type 2 are commuted to each other and so the commuting conjugacy class graph has a subgraph isomorphic to $K_{mp(p-1)}$. We now determine the relationship between conjugacy classes of this and other types. To do this, we assume that $a^i$ is a representatives of a conjugacy class of Type 2 and $a^ub^vz$ is an arbitrary element of $G$ such that $(a^i)(a^ub^vz) = (a^ub^vz)(a^i)$. So  $a^ib^v = b^va^i$ and since $p\nmid i$, by Lemma \ref{lem3.2.5}(a) the latter is true if and only if $v=0$. Therefore, $a^i$ commutes with all elements in the form $a^uz$. This shows that the conjugacy classes of Type 2 are commuted only with the conjugacy classes of Type 1.

\item Suppose $b^{t_1p}$ and $b^{t_2p}$ are the representative of two conjugacy classes of Type 3. It is clear that $b^{t_1p}b^{t_2p}=b^{t_2p}b^{t_1p}$ and so all conjugacy classes of Type 3 are commuting together in the graph. Hence by Table 1, the graph has a subgraph isomorphic to $K_{n(p-1)}$. To determine the relationship between conjugacy classes of this and other types, we choose $b^{tp}$ to be a representative of a conjugacy class of Type 3  and $a^ub^vz$ as an arbitrary element of $G$ such that $(b^{tp})(a^ub^vz) = (a^ub^vz)(b^{tp})$. Thus, $b^{tp}a^u = a^ub^{tp}$ and by Lemma \ref{lem3.2.5}(a), the later is true if and only if $u \equiv 0 \ (mod \ p)$ or $u=sp$ in which $0\leq s\leq p-1$. Therefore, $b^{tp}$ can be commuted with all elements in the form $a^{sp}b^vz$. Hence, the following are satisfied:
\begin{enumerate}
\item If $s=0$, then the conjugacy classes of Types 3 and 4 are commuting together.

\item If $s\neq 0$ and $v=0$, then the conjugacy classes of Types 3 and 1 are commuting to each other.

\item If $v,s\neq 0$ and $p\mid v$, then the conjugacy classes of Types 3 and 5 are commuting together.

\item If $s\neq 0$ and $p\nmid v$, then the conjugacy classes of Types 3 and 7 are commuting to each other.
\end{enumerate}

\item Suppose $b^{j_1}$ and $b^{j_2}$ are the representatives of two conjugacy classes of Type 4.
It is obvious that $b^{j_1}b^{j_2}=b^{j_2}b^{j_1}$ which shows that all conjugacy classes of Type 4 are commuted to each other. So, the commuting conjugacy class graph has a subgraph isomorphic to $K_{mp(p-1)}$. To determine the relationship between this and other types in Table 1, we assume that $b^j$ is a representative of a conjugacy class of type 4  and $a^ub^vz$ is an arbitrary element of $G$ such that $(b^j)(a^ub^vz) = (a^ub^vz)(b^j)$. Thus,  $b^ja^u = a^ub^j$ and since $p\nmid j$, by Lemma \ref{lem3.2.5}(a), the latter is true if and only if $u=0$. Therefore, $b^j$ commutes with all elements in the form $b^vz$ and hence the conjugacy classes of Type 4 can be commuted only with the conjugacy classes of Type 3.

\item Suppose $a^{s_1p}b^{t_1p}$ and $a^{s_2p}b^{t_2p}$ are the representatives of two conjugacy classes of Type 5. By Equation \ref{e45}, $(a^{s_1p}b^{t_1p})(a^{s_2p}b^{t_1p})=(a^{s_2p}b^{t_2p})(a^{s_1p}b^{t_2p})$ and so all conjugacy classes of Type 5 are commuted to each other. This gives us the complete graph $K_{n(p-1)^2}$ as a subgraph of  the commuting conjugacy class graph. Hence, it is enough to determine the relationship between this and conjugacy classes of Types 6, 7 and 8. Suppose $a^{sp}b^{tp}$ and $a^ub^{vp}$ are the representatives of the conjugacy classes of Types 5 and 6, respectively, such that they are commuting together. By Lemma \ref{lem2.4}, $a^{sp}b^{tp}=(a^{ip}b^p)^{n_t}z_1$ and $a^ub^{vp}=(ab^{jp})^{n_u}z_2$ and so by Lemma \ref{lem3.2}, $C_G(a^{sp}b^{tp})=C_G(a^{ip}b^p)$ and $C_G(a^ub^{vp})=C_G(ab^{jp})$. It is now easy to prove that  $(a^{sp}b^{tp})(a^ub^{vp})=(a^ub^{vp})(a^{sp}b^{tp})$ if and only if $(a^{ip}b^p)(ab^{jp})=(ab^{jp})(a^{ip}b^p)$ if and only if $b^pa=ab^p$ which is a contradiction with Lemma \ref{lem3.2.5}(a). Therefore, a conjugacy class of Type 5 is not adjacent with a conjugacy class of Type 6. Similarly, a conjugacy class of Type 5 is not adjacent to another one of Type 7. We now assume that $a^{sp}isb^{tp}$ and $a^ub^v$ are representatives of the conjugacy classes of Types 5 and 8, respectively. By Lemma \ref{lem2.4}, $a^{sp}b^{tp}=(a^{ip}b^p)^{n_t}z_1$ and $a^ub^v=(a^jb)^{n_v}z_2$. On the other hand, By Lemma \ref{lem3.2}, $C_G(a^{sp}b^{tp})=C_G(a^{ip}b^p)$ and $C_G(a^ub^v)=C_G(a^jb)$ and by Equation \ref{e12-5}, we can see $(a^{sp}b^{tp})(a^ub^v)=(a^ub^v)(a^{sp}b^{tp})$ if and only if $(a^{ip}b^p)(a^jb)=(a^jb)(a^{ip}b^p)$ if and only if  $j\equiv i (mod \ p)$. Since $1\leq i\leq p-1$, the conjugacy classes of Types 5 and 8 can be divided into $p-1$ parts. Moreover, this process constructs a complete subgraph of size $n(p-1)$ in the commuting conjugacy class graph.

\item Suppose  $a^{i_1}b^{t_1p}$ and $a^{i_2}b^{t_2p}$ are the representative of two conjugacy classes of Types 6 such that they are commute to each other. By Lemma \ref{lem2.4}, $a^{i_1}b^{t_1p}=(ab^{s_1p})^{n_{i_1}}z_1$ and $a^{i_2}b^{t_2p}=(ab^{s_2p})^{n_{i_2}}z_2$. On the other hand, by Lemma (\ref{lem3.2}), $C_G(a^{i_1}b^{t_1p})=C_G(ab^{s_1p})$ and $C_G(a^{i_2}b^{t_2p})=C_G(ab^{s_2p})$. Therefore, $(a^{i_1}b^{t_1p})(a^{i_2}b^{t_2p})=(a^{i_2}b^{t_2p})(a^{i_1}b^{t_1p})$ if and only if $(ab^{s_1p})(ab^{s_2p})=(ab^{s_2p})(ab^{s_1p})$ if and only if $ab^{(s_2-s_1)p}=b^{(s_2-s_1)p}a$. By Lemma  \ref{lem3.2.5}(a), the last equality is satisfied if and only if $s_1=s_2$. This proves that in the commuting conjugacy class  graph, two conjugacy classes of Type 6  are adjacent, when the centralizers of their representatives is equal. By Theorem \ref{thm3.5}, the number of centralizer of Types 6 is $p-1$ and so the conjugacy classes of this type can be divided into $p-1$ parts. On the other hand, by Table \ref{j3.1},  the number of conjugacy classes of Type 6 is $mp(p-1)^2$. Hence, each part of Type 6 is a clique of size $mp(p-1)$. Now it is enough to determine the relationship between conjugacy classes of Type 6 with other conjugacy classes of Types 7 and 8.  To do this, we assume that $a^ib^{tp}$ and $a^{sp}b^j$  are the representatives of conjugacy classes of Types 6 and 7, respectively, such that they are commute to each other. By Equation \ref{e45}, $(a^ib^{tp})(a^{sp}b^j)=(a^{sp}b^j)(a^ib^{tp})$ if and only if $a^ib^j = b^ja^i$. Since $p\nmid i,j$, Lemma \ref{lem3.2.5}(a) implies that the last equality  cannot be occurred. This means that a conjugacy class of Type 6 is not adjacent with another one of Type 7 in the commuting conjugacy class graph. Next, we assume that $a^ib^{tp}$ and $a^ub^v$ are the representatives of two conjugacy classes of Types 6 and 8, respectively and they are commute to each other. By Lemma \ref{lem2.4}, $a^ib^{tp}=(ab^{sp})^{n_i}z_1$ and $a^ub^v=(ab^j)^{n_u}z_2$ and  by Lemma \ref{lem3.2}, $C_G(a^ib^{tp})=C_G(ab^{sp})$ and $C_G(a^ub^v)=C_G(ab^j)$. Therefore,  $(a^ib^{tp})(a^ub^v)=(a^ub^v)(a^ib^{tp})$ if an only if $(ab^{sp})(ab^j)=(ab^j)(ab^{sp})$ if and only if $b^{sp-j}a=ab^{sp-j}$. We now apply Lemma \ref{lem3.2.5}(a) to prove that the last equality is satisfied  if and only if $p^2\mid (sp-j)$ or $p\mid j$ which is a contradiction. Since $p\nmid j$, the conjugacy classes of Type  6 is not adjacent with another one of Type 8 in the commuting conjugacy class  graph.

\item Suppose that $a^{s_1p}b^{j_1}$ and $a^{s_2p}b^{j_2}$ are representatives of two conjugacy classes of Type 7  such that they are commuting to each other. By Lemma \ref{lem2.4}, $a^{s_1p}b^{j_1}=(a^{t_1p}b)^{n_{j_1}}z_1$ and $a^{s_2p}b^{j_2}=(a^{t_2p}b)^{n_{j_2}}z_2$ and by Lemma \ref{lem3.2}, $C_G(a^{s_1p}b^{j_1})=C_G(a^{t_1p}b)$ and $C_G(a^{s_2p}b^{j_2})=C_G(a^{t_2p}b)$. Therefore,  $(a^{s_1p}b^{j_1})(a^{s_2p}b^{j_2})=(a^{s_2p}b^{j_2})(a^{s_1p}b^{j_1})$ if and only if $(a^{t_1p}b)(a^{t_2p}b)=(a^{t_2p}b)(a^{t_1p}b)$ if and only if $ba^{(t_2-t_1)p}=a^{(t_2-t_1)p}b$. By Lemma  \ref{lem3.2.5}(a), the last equality is satisfied if and only if $t_1=t_2$. This means that two conjugacy classes of Types 7 are adjacent in the commuting conjugacy class graph, when their centralizers is equal. Also, by Theorem \ref{thm3.5}, the number of centralizers of Type 7 is $p-1$. So, the conjugacy classes of this type can be divided into $p-1$ parts. Now by Table \ref{j3.1},  the number of conjugacy classes of Type 7 is equal to  $mp(p-1)^2$ and the  graph structure of each part  of Type 7 gives a complete subgraph of size $mp(p-1)$. It is enough to find the relationship between this and other conjugacy classes of Type 8. To do this, we assume that $a^{sp}b^j$ and $a^ub^v$ are representatives of two conjugacy classes of Types 7 and 8, respectively such that they are commute to each other. By Lemma \ref{lem2.4}, $a^{sp}b^j=(a^{tp}b)^{n_j}z_1$ and $a^ub^v=(a^ib)^{n_v}z_2$. On the other hand, by Lemma \ref{lem3.2},  $C_G(a^{sp}b^j)=C_G(a^{tp}b)$ and $C_G(a^ub^v)=C_G(a^ib)$. Therefore, $(a^{sp}b^j)(a^ub^v)=(a^ub^v)(a^{sp}b^j)$ if and only if $(a^{tp}b)(a^ib)=(a^ib)(a^{tp}b)$ if and only if $a^{tp-i}b=ba^{tp-i}$. By Lemma \ref{lem3.2.5}(a), the last equality is satisfied if and only if $p^2\mid (tp-i)$ or $p\mid i$ which is a contradiction. Thus, the conjugacy classes of Type 7 is not adjacent with a conjugacy class of Type 8 in the commuting conjugacy class graph.

\item Suppose  $a^{i_1}b^{j_1}$ and $a^{i_2}b^{j_2}$ are the representatives of two conjugacy classes of Type 8 such that they are commuting together. By Lemma \ref{lem2.4}, $a^{i_1}b^{j_1}=(a^sb)^{n_{j_1}}z_1$ and $a^{i_2}b^{j_2}=(a^tb)^{n_{j_2}}z_2$ and by Lemma \ref{lem3.2}, $C_G(a^{i_1}b^{j_1})=C_G(a^sb)$ and $C_G(a^{i_2}b^{j_2})=C_G(a^tb)$. Therefore, $(a^{i_1}b^{j_1})(a^{i_2}b^{j_2})=(a^{i_2}b^{j_2})(a^{i_1}b^{j_1})$ if and only if $(a^sb)(a^tb)=(a^tb)(a^sb)$ if and only if $a^{s-t}b=ba^{s-t}$. By Lemma  \ref{lem3.2.5}(a), the last equality is satisfied if and only if $s=t$. This means that two conjugacy classes of Type 8  are adjacent in the commuting conjugacy class graph, when their centralizers is equal. Also, by Theorem \ref{thm3.5}, the number of centralizers of Type 8 is $p(p-1)$ and so the conjugacy classes of this type can be divided into $p(p-1)$ parts. On the other hand, by Table \ref{j3.1},  the number of conjugacy classes of Type 8 is  $mp^2(p-1)^2$ which proves that each part of Type 8 gives a clique of size $mp(p-1)$ in the commuting conjugacy class graph. By our discussion in Case 5, all conjugacy classes of Type 8 are adjacent only with conjugacy classes of Type 5. Since the conjugacy classes of Types 5  and  8 can be divided into $p-1$ and $p(p-1)$ parts, every part of Type 5  is adjacent with $p$ parts of Type 8.
\end{enumerate}

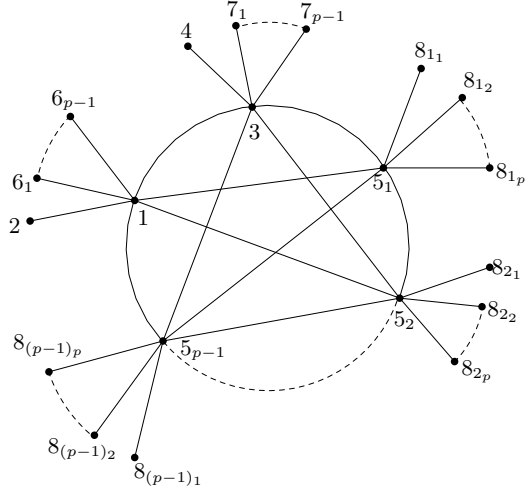
\begin{figure}[h]
\begin{picture}(400,300)
\scalebox{.95}{
\definecolor{qqqqff}{rgb}{0,0,0}
\begin{tikzpicture}[line cap=round,line join=round,>=triangle 45,x=1.2cm,y=1.2cm]
\clip(-5.66,-4.54) rectangle (15.68,4.94);
\draw [shift={(0.55,1.46)}] plot[domain=-0.37:3.87,variable=\t]({1*1.65*cos(\t r)+0*1.65*sin(\t r)},{0*1.65*cos(\t r)+1*1.65*sin(\t r)});
\draw [shift={(0.56,1.38)},dash pattern=on 2pt off 2pt]  plot[domain=3.84:5.96,variable=\t]({1*1.6*cos(\t r)+0*1.6*sin(\t r)},{0*1.6*cos(\t r)+1*1.6*sin(\t r)});
\draw [shift={(0.57,2.81)},dash pattern=on 2pt off 2pt]  plot[domain=1.22:1.88,variable=\t]({1*1.27*cos(\t r)+0*1.27*sin(\t r)},{0*1.27*cos(\t r)+1*1.27*sin(\t r)});
\draw [shift={(1.69,2.22)},dash pattern=on 2pt off 2pt]  plot[domain=0.11:0.71,variable=\t]({1*1.46*cos(\t r)+0*1.46*sin(\t r)},{0*1.46*cos(\t r)+1*1.46*sin(\t r)});
\draw [shift={(1.84,0.98)},dash pattern=on 2pt off 2pt]  plot[domain=5.51:6.1,variable=\t]({1*1.24*cos(\t r)+0*1.24*sin(\t r)},{0*1.24*cos(\t r)+1*1.24*sin(\t r)});
\draw [shift={(-0.42,0.53)},dash pattern=on 2pt off 2pt]  plot[domain=3.46:4.02,variable=\t]({1*1.67*cos(\t r)+0*1.67*sin(\t r)},{0*1.67*cos(\t r)+1*1.67*sin(\t r)});
\draw [shift={(-0.62,1.93)},dash pattern=on 2pt off 2pt]  plot[domain=2.39:2.93,variable=\t]({1*1.55*cos(\t r)+0*1.55*sin(\t r)},{0*1.55*cos(\t r)+1*1.55*sin(\t r)});
\draw (-1,2)-- (-2.22,1.76);
\draw (-2.14,2.26)-- (-1,2);
\draw (-1.75,2.98)-- (-1,2);
\draw (-0.38,3.8)-- (0.37,3.09);
\draw (0.18,4.04)-- (0.37,3.09);
\draw (1,4)-- (0.37,3.09);
\draw (2.34,3.54)-- (1.9,2.38);
\draw (2.82,3.2)-- (1.9,2.38);
\draw (3.14,2.38)-- (1.9,2.38);
\draw (3.14,1.22)-- (2.09,0.86);
\draw (3.05,0.76)-- (2.09,0.86);
\draw (2.73,0.12)-- (2.09,0.86);
\draw (-1,-1)-- (-0.67,0.36);
\draw (-1.47,-0.74)-- (-0.67,0.36);
\draw (-2,0)-- (-0.67,0.36);
\draw (1.9,2.38)-- (-1,2);
\draw (-1,2)-- (2.09,0.86);
\draw (2.09,0.86)-- (-0.67,0.36);
\draw (1.9,2.38)-- (-0.67,0.36);
\draw (0.37,3.09)-- (-0.67,0.36);
\draw (0.37,3.09)-- (2.09,0.86);
\begin{scriptsize}
\fill [color=qqqqff] (-1,2) circle (1.5pt);
\draw[color=qqqqff] (-0.9,1.8) node {$1$};
\fill [color=qqqqff] (-2.22,1.76) circle (1.5pt);
\draw[color=qqqqff] (-2.4,1.7) node {$2$};
\fill [color=qqqqff] (-2.14,2.26) circle (1.5pt);
\draw[color=qqqqff] (-2.3,2.2) node {$6_1$};
\fill [color=qqqqff] (-1.75,2.98) circle (1.5pt);
\draw[color=qqqqff] (-1.7,3.17) node {$6_{p-1}$};
\fill [color=qqqqff] (0.37,3.09) circle (1.5pt);
\draw[color=qqqqff] (0.4,2.8) node {$3$};
\fill [color=qqqqff] (-0.38,3.8) circle (1.5pt);
\draw[color=qqqqff] (-0.4,3.98) node {$4$};
\fill [color=qqqqff] (0.18,4.04) circle (1.5pt);
\draw[color=qqqqff] (0.2,4.22) node {$7_1$};
\fill [color=qqqqff] (1,4) circle (1.5pt);
\draw[color=qqqqff] (1.2,4.18) node {$7_{p-1}$};
\fill [color=qqqqff] (1.9,2.38) circle (1.5pt);
\draw[color=qqqqff] (1.9,2.2) node {$5_1$};
\fill [color=qqqqff] (2.34,3.54) circle (1.5pt);
\draw[color=qqqqff] (2.45,3.73) node {$8_{1_1}$};
\fill [color=qqqqff] (2.82,3.2) circle (1.5pt);
\draw[color=qqqqff] (3.02,3.38) node {$8_{1_2}$};
\fill [color=qqqqff] (3.14,2.38) circle (1.5pt);
\draw[color=qqqqff] (3.4,2.3) node {$8_{1_p}$};
\fill [color=qqqqff] (2.09,0.86) circle (1.5pt);
\draw[color=qqqqff] (2.15,.6) node {$5_2$};
\fill [color=qqqqff] (3.14,1.22) circle (1.5pt);
\draw[color=qqqqff] (3.35,1.2) node {$8_{2_1}$};
\fill [color=qqqqff] (3.05,0.76) circle (1.5pt);
\draw[color=qqqqff] (3.3,0.7) node {$8_{2_2}$};
\fill [color=qqqqff] (2.73,0.12) circle (1.5pt);
\draw[color=qqqqff] (3.,0.0) node {$8_{2_p}$};
\fill [color=qqqqff] (-0.67,0.36) circle (1.5pt);
\draw[color=qqqqff] (-0.2,0.25) node {$5_{p-1}$};
\fill [color=qqqqff] (-1,-1) circle (1.5pt);
\draw[color=qqqqff] (-0.6,-1.2) node {$8_{(p-1)_1}$};
\fill [color=qqqqff] (-1.47,-0.74) circle (1.5pt);
\draw[color=qqqqff] (-1.6,-0.9) node {$8_{(p-1)_2}$};
\fill [color=qqqqff] (-2,0) circle (1.5pt);
\draw[color=qqqqff] (-2.,0.3) node {$8_{(p-1)_p}$};
\end{scriptsize}
\end{tikzpicture}
}
\end{picture}
\vspace*{-3cm}
\caption{The Commuting Conjugacy Class Graph of $G$ with $\frac{G}{Z(G)}\cong \mathbb{Z}_{p^2} \times \mathbb{Z}_{p^2}$.}
\label{fig4}
\end{figure}

By above discussion in eight parts, the commuting conjugacy class graph of $G$  is  a connected graph with $n(p^2-1)+m(p^4-p^2)$ vertices. Suppose $\mathcal{M}_1$ is a graph depicted in Figure \ref{fig3}. Therefore, the commuting conjugacy class graph of $G$ is a $\mathcal{M}_1$-join of graphs as follows:
\[\Gamma(G) = \mathcal{M}_1[\overbrace{K_{n(p-1)},\cdots,K_{n(p-1)}}^{p+1},\overbrace{K_{m(p^2-p)},\cdots, K_{m(p^2-p)}}^{p^2+p}],\]
where $n=\frac{|Z|}{p}$ and $m=\frac{|Z|}{p^2}$. The graph $\Gamma(G)$ is depicted in Figure \ref{fig4}.
\end{proof}

In the following theorem, we study the commuting conjugacy class graph of $G$ with this condition that $\frac{G}{Z(G)}$ $\cong$ $ \mathbb{Z}_{p^2}\rtimes \mathbb{Z}_{p^2}$, when the semidirect product is non-abelian. We recall that by Theorem \ref{th2.0.0}, $p$ is an odd prime number.

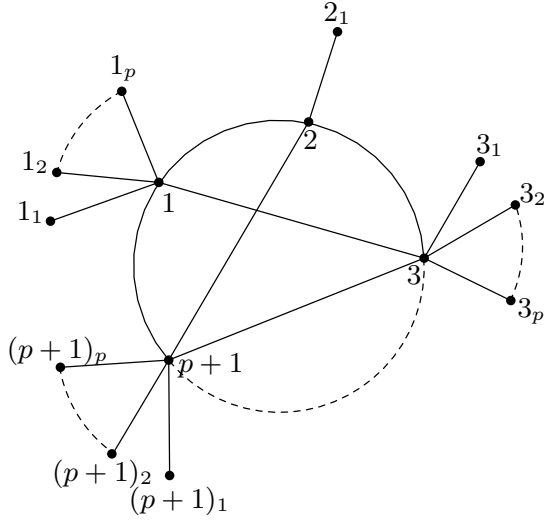
\begin{figure}[h]
\begin{picture}(350,200)
\definecolor{qqqqff}{rgb}{0,0,0}
\scalebox{1.2}{
\begin{tikzpicture}[line cap=round,line join=round,>=triangle 45,x=1.0cm,y=1.0cm]
\clip(-5.66,-4.54) rectangle (15.68,4.94);
\draw [shift={(0.55,1.41)}] plot[domain=0.05:3.86,variable=\t]({1*1.61*cos(\t r)+0*1.61*sin(\t r)},{0*1.61*cos(\t r)+1*1.61*sin(\t r)});
\draw [shift={(0.56,1.38)},dash pattern=on 2pt off 2pt]  plot[domain=-2.44:0.07,variable=\t]({1*1.6*cos(\t r)+0*1.6*sin(\t r)},{0*1.6*cos(\t r)+1*1.6*sin(\t r)});
\draw [shift={(1.86,1.62)},dash pattern=on 2pt off 2pt]  plot[domain=-0.45:0.33,variable=\t]({1*1.39*cos(\t r)+0*1.39*sin(\t r)},{0*1.39*cos(\t r)+1*1.39*sin(\t r)});
\draw [shift={(-0.42,0.53)},dash pattern=on 2pt off 2pt]  plot[domain=3.31:4.09,variable=\t]({1*1.48*cos(\t r)+0*1.48*sin(\t r)},{0*1.48*cos(\t r)+1*1.48*sin(\t r)});
\draw [shift={(-0.51,2.06)},dash pattern=on 2pt off 2pt]  plot[domain=2.06:2.88,variable=\t]({1*1.45*cos(\t r)+0*1.45*sin(\t r)},{0*1.45*cos(\t r)+1*1.45*sin(\t r)});
\draw (-0.78,2.33)-- (-1.98,1.9);
\draw (-1.91,2.44)-- (-0.78,2.33);
\draw (-1.19,3.34)-- (-0.78,2.33);
\draw (1.2,4.)-- (0.88,3);
\draw (-0.66,-0.92)-- (-0.67,0.36);
\draw (-1.29,-0.68)-- (-0.67,0.36);
\draw (-1.87,0.28)-- (-0.67,0.36);
\draw (0.88,3)-- (-0.67,0.36);
\draw (3.12,1.02)-- (2.16,1.49);
\draw (3.17,2.08)-- (2.16,1.49);
\draw (2.78,2.56)-- (2.16,1.49);
\draw (-0.78,2.33)-- (2.16,1.49);
\draw (2.16,1.49)-- (-0.67,0.36);
\begin{scriptsize}
\fill [color=qqqqff] (-0.78,2.33) circle (1.5pt);
\draw[color=qqqqff] (-0.67,2.1) node {$1$};
\fill [color=qqqqff] (-1.19,3.34) circle (1.5pt);
\draw[color=qqqqff] (-1.17,3.6) node {$1_p$};
\fill [color=qqqqff] (-1.91,2.44) circle (1.5pt);
\draw[color=qqqqff] (-2.15,2.55) node {$1_2$};
\fill [color=qqqqff] (-1.98,1.9) circle (1.5pt);
\draw[color=qqqqff] (-2.2,2.0) node {$1_1$};
\fill [color=qqqqff] (0.88,3) circle (1.5pt);
\draw[color=qqqqff] (0.9,2.8) node {$2$};
\fill [color=qqqqff] (1.2,4.) circle (1.5pt);
\draw[color=qqqqff] (1.2,4.2) node {$2_1$};
\fill [color=qqqqff] (2.16,1.49) circle (1.5pt);
\draw[color=qqqqff] (2.05,1.3) node {$3$};
\fill [color=qqqqff] (3.12,1.02) circle (1.5pt);
\draw[color=qqqqff] (3.35,.9) node {$3_p$};
\fill [color=qqqqff] (3.17,2.08) circle (1.5pt);
\draw[color=qqqqff] (3.37,2.2) node {$3_2$};
\fill [color=qqqqff] (2.78,2.56) circle (1.5pt);
\draw[color=qqqqff] (2.89,2.74) node {$3_1$};
\fill [color=qqqqff] (-0.67,0.36) circle (1.5pt);
\draw[color=qqqqff] (-0.2,0.3) node {$p+1$};
\fill [color=qqqqff] (-1.87,0.28) circle (1.5pt);
\draw[color=qqqqff] (-1.9,0.46) node {$(p+1)_p$};
\fill [color=qqqqff] (-1.3,-0.68) circle (1.5pt);
\draw[color=qqqqff] (-1.4,-0.9) node {$(p+1)_2$};
\fill [color=qqqqff] (-0.66,-0.92) circle (1.5pt);
\draw[color=qqqqff] (-0.55,-1.2) node {$(p+1)_1$};
\end{scriptsize}
\end{tikzpicture}
}
\end{picture}
\vspace*{-3cm}
\caption{The Graph $\mathcal{M}_2$.}\label{fig5}
\end{figure}

\begin{thm}
Suppose $p$ is an odd prime number and $G$ is a group with center $Z$ such that $\frac{G}{Z}\cong \mathbb{Z}_{p^2}\rtimes \mathbb{Z}_{p^2}$ and $\frac{G}{Z}$ is non-abelian. Then the commuting conjugacy classes of $G$ is as follows:
\[\mathcal{M}_2[\overbrace{K_{n(p-1)},\cdots,K_{n(p-1)}}^{p^2+p+1},K_{np(p-1)}],\]
where $\mathcal{M}_2$ is depicted in Figure \ref{fig5}. Here, $n = \frac{|Z|}{p}$.
\end{thm}

\vskip 20mm

\begin{proof}

We first apply  Equation \ref{e1} to deduce that
\begin{equation}\label{e15}
b^ja^iZ = a^{ijp+i}b^jZ
\end{equation}
which implies that
\begin{eqnarray}
b^jZ a^iZ & = & (a^{p}Z)^{ij}a^{i}Zb^jZ   \label{c2},\\
a^iZb^jZ  & = & (a^{p}Z)^{-ij}b^jZa^{i}Z  \label{c3}.
\end{eqnarray}

Choose elements $x$ and $g$ in $G\setminus Z$. By the proof of Theorem \ref{thm3.5}, $x=a^ib^jz_1$ and $g=a^ub^vz_2$, where $z_1,z_2\in Z$, $0\leq i, j, u, v\leq p^2-1$ and $i , j$, as well as $u , v$, are not simultaneously zero. In Equation \ref{e15}, we put  $i=p$ and $j=1$. Since $a^{p^2}\in Z$, $bZa^pZ=a^pZbZ$. Thus  $a^pZ\in Z\biggl(\frac{G}{Z}\biggr)$ and by Equations  \ref{c2}  and \ref{c3},

\begin{eqnarray*}
g^{-1}xgZ & = & g^{-1}ZxZgZ \\
               & = & (a^ub^v)^{-1}Z(a^ib^j)Z(a^ub^v)Z \\
               & = & (a^ub^v)^{-1}Za^iZb^jZa^uZb^vZ \\
               & = & (a^ub^v)^{-1}Za^iZ(a^{p}Z)^{uj}a^uZb^jZb^vZ   \\
               & = & (a^{p}Z)^{uj}(a^ub^v)^{-1}Za^uZa^iZb^vZb^jZ   \\
               & = & (a^{p}Z)^{uj}(a^ub^v)^{-1}Za^uZ(a^{p}Z)^{-vi}b^vZa^iZb^jZ \\
               & = & (a^{p}Z)^{uj}(a^{p}Z)^{-vi}(a^ub^v)^{-1}Za^uZb^vZa^iZb^jZ \\
               & = & (a^{p}Z)^{uj-vi}(a^ub^v)^{-1}Z(a^ub^v)Z(a^ib^j)Z \\
               & = & (a^{(uj-vi)p}Z)(a^ib^j)Z \\
               & = & (a^{(uj-vi)p}a^ib^j)Z.
\end{eqnarray*}
So, there exists $z\in Z$ such that
\begin{equation}\label{c4}
(a^ub^v)^{-1}(a^ib^j)(a^ub^v)=a^{(uj-vi)p}a^ib^jz.
\end{equation}
Suppose that $x = a^ib^jz$ is an arbitrary element of $G$ such that $0 \leq i,j \leq p^2 -1$ and $i,j$ are not simultaneously zero.  To compute the number of non-central conjugacy classes of $G$ and obtain their structure, we will consider  three separate cases that all together can be divided into eight subcases:

\begin{enumerate}[a)]
\item \textit{$i\neq 0$ and $j=0$.} This case is separated into two subcases that $p \mid i$ and $p \nmid i$ as follows:

\begin{enumerate}[1)]
\item  $p\mid i$. Suppose $i=sp$, where $1\leq s\leq p-1$. By Equation \ref{c4}  and this fact that $a^{p^2}\in Z$,  $(a^ub^v)^{-1}(a^{sp})(a^ub^v)=a^{-vsp^2}a^{sp}z_1=a^{sp}z_2$ and so $(a^{sp})^G=\{a^{sp}z\mid z\in Z\}$. On the other hand, by Theorem \ref{thm3.5},  $|C_G(a^{sp})|=p^3|Z|$ and $|(a^{sp})^G|=p$. Thus, $(a^{sp})^G=\{a^{sp}z_1,\cdots,a^{sp}z_p\}=a^{sp}H$ in which $H=\{z_1,\cdots,z_p\}\subseteq Z$. Suppose $|Z|>p$ and choose the elements $z_r\in Z\setminus H$ and an element in the form  $a^{sp}z_r$. It is easy to see that $a^{sp}z_r\not\in (a^{sp})^G$ and hence $(a^{sp}z_r)^G\neq (a^{sp})^G$. Furthermore, $C(a^{sp}z_r)=C(a^{sp})$ and so $|(a^{sp}z_r)^G|=|(a^{sp})^G|=p$. Also, $(a^{sp}z_{r_1})^G=(a^{sp})^Gz_{r_1}=a^{sp}Hz_{r_1}$. Since $(a^{sp})^G \cap (a^{sp}z_{r_1})^G = \emptyset$, $H\cap Hz_{r_1}=\emptyset$ and $H\cup Hz_{r_1}\subseteq Z$. We now choose the element $z_{r_2}\in Z\setminus(H\cup Hz_{r_1})$ and continue this process to prove that $p\mid |Z|$.  If  $n=\frac{|Z|}{p}$ then there are $n$ distinct conjugacy classes in the form of $(a^{sp}z_{r_1})^G,(a^{sp}z_{r_2})^G\cdots (a^{sp}z_{r_n})^G$. Since $1\leq s\leq p-1$, there are $n(p-1)$ distinct conjugacy classes with  $n=\frac{|Z|}{p}$ and each conjugacy class has $p$ elements.

\item $p\nmid i$. By Equation \ref{c4},
\begin{equation}\label{c5}
(a^ub^v)^{-1}(a^{i})(a^ub^v)=a^{-vip}a^{i}z.
\end{equation}
By Theorem 3.2,  $|C_G(a^{i})|=p^2|Z|$ and so $|(a^{i})^G|=p^2$. Note that $0\leq v\leq p^2-1$ and so there are $v'$ and $k$ such that $-vi=v'p+k$ and $0\leq k\leq p-1$. By Equation \ref{c5}  and the fact that $a^{p^2}\in Z$,
\begin{equation}\label{c6}
(a^{i})^G=\{a^{kp}a^iz_{(k,l)}\mid 0\leq k,l\leq p-1,z_{(k,l)}\in Z\}.
\end{equation}

Set $A_i=\{a^{kp}a^iz\mid 0\leq k\leq p-1, z\in Z\}$ and $H=\{z_{(0,l)}\mid 0\leq l\leq p-1\}$. Obviously, $|A_i|=p|Z|$, $(a^i)^G\subseteq A_i$ and $H\subseteq Z$. If $|Z|>p$, then  we choose an element $z_r\in Z\setminus H$ and an element  $a^{i}z_r \in A_i$. It can be easily seen that   $a^{i}z_r\not\in (a^{i})^G$ and $(a^{i}z_r)^G\neq (a^{i})^G$. Since $C(a^{i}z_r)=C(a^{i})$,   $|(a^{i}z_r)^G|=|(a^{i})^G|=p^2$ and
$(a^iz_{r_1})^G=(a^i)^Gz_{r_1}$. On the other hand, since $(a^i)^G\cap(a^iz_{r_1})^G=\emptyset$, $H\cap Hz_{r_1}=\emptyset$ and $H\cup Hz_{r_1}\subseteq Z$. We now choose an element $z_{r_2}\in Z\setminus(H\cup Hz_{r_1})$ and continue this process to prove that $(a^{i}z_{r_2})^G\subseteq A_i$. This shows that $p^2\mid |A_i|$. Put $n=\frac{|A_i|}{p^2}=\frac{|Z|}{p}$. Then, for each $i$ there are $n$ distinct conjugacy classes as $(a^{i}z_{r_1})^G,\cdots,(a^{i}z_{r_n})^G$. Again since $1\leq i\leq p^2-1$ and $p\nmid i$,  there are non-negative integers $k'$ and $i'$ such that  $i=k'p+i'$ and $1\leq i'\leq p-1$. By Equation  \ref{c6},
\begin{eqnarray*}			
(a^{i})^G & = & (a^{k'p+i'})^G \\
               & = & \{a^{kp}a^{k'p+i'}z_{(k,l)}\mid 0\leq k,l\leq p-1,z_{(k,l)}\in Z\} \\
               & = & \{a^{(k+k')p}a^{i'}z_{(k,l)}\mid 0\leq k,l\leq p-1,z_{(k,l)}\in Z\}.
\end{eqnarray*}
Therefore, there exists $z_{i'}\in Z$ such that $(a^{i})^G=(a^{i'}z_{i'})^G$. This proves that  the number of distinct conjugacy classes   is equal to $n(p-1)$ in which $n=\frac{|Z|}{p}$ and each class has $p^2$ elements.
\end{enumerate}

\item \textit{$i = 0$ and $j \neq 0$.} Again there are two cases that $p \mid j$ and $p \nmid j$.
\begin{enumerate}[1)]
\setcounter{enumii}{2}
\item $p\mid j$. Suppose $j=tp$, where $1\leq t\leq p-1$. By Equation \ref{c4} and the fact that $a^{p^2}$ is an element of $Z$, we can conclude that $(a^ub^v)^{-1}(b^{tp})(a^ub^v)=a^{utp^2}b^{tp}z_1=b^{tp}z_2$ and so $(b^{tp})^G=\{b^{tp}z\mid z\in Z\}$. By  Theorem \ref{thm3.5}, $|C_G(b^{tp})|=p^3|Z|$ and hence $|(b^{tp})^G|=p$. Therefore, similar to the Case 1, there are $n$ distinct conjugacy classes $(b^{tp}z_{r_1})^G,\cdots,(b^{tp}z_{r_n})^G$, where $n=\frac{|Z|}{p}$. So, the number of distinct conjugacy classes is equal to $n(p-1)$  and each class has $p$ elements.

\item $p\nmid j$. Then by Equation \ref{c4}, $(a^ub^v)^{-1}(b^{j})(a^ub^v)=a^{ujp}b^{j}z$.
By Theorem \ref{thm3.5},  $|C_G(b^{j})|=p^2|Z|$ and so $|(b^{j})^G|=p^2$. A similar argument as in the Case 2 shows that $(b^{j})^G=\{a^{kp}b^jz_{(k,l)}\mid 0\leq k,l\leq p-1,z_{(k,l)}\in Z\}$. Thus, there are $n$ distinct conjugacy classes $(b^{j}z_{r_1})^G,\cdots,(b^{j}z_{r_n})^G$, when $j$ is a fixed number and $n=\frac{|Z|}{p}$. Since $1\leq j\leq p^2-1$ and $p\nmid j$, the number of distinct conjugacy classes is equal to $n(p^2-p)$  and each class has $p^2$ elements.
\end{enumerate}

\item \textit{$i\neq 0$ and $j \neq 0$.} We have four subcases as follows:
\begin{enumerate}[1)]
\setcounter{enumii}{4}
\item \textit{$p\mid i$ and $p\mid j$}. Suppose $i=sp$ and $j=tp$ such that $1\leq s,t\leq p-1$. By Equation \ref{c4} and the fact that $a^{p^2}\in Z$, $(a^ub^v)^{-1}(a^{sp}b^{tp})(a^ub^v)$ $=$ $a^{(utp-vsp)p}a^{sp}b^{tp}z_1$ $=$ $a^{(ut-vs)p^2}a^{sp}b^{tp}z_1$ $=$ $a^{sp}b^{tp}z_2.$ Thus, $(a^{sp}b^{tp})^G=\{a^{sp}b^{tp}z\mid z\in Z\}.$ Apply Theorem \ref{thm3.5} to deduce that $|C_G(a^{sp}b^{tp})|=p^3|Z|$ which implies that  $|(a^{sp}b^{tp})^G|=p$. Similar to the Case 1, for fixed non-negative integers $s$ and $t$, there are $n$ distinct conjugacy classes as  $(a^{sp}b^{tp}z_{r_1})^G$, $\ldots$, $(a^{sp}b^{tp}z_{r_n})^G$ in which $n = \frac{|Z|}{p}$. By above discussion, we result that the number of the distinct conjugacy classes is equal to $n(p-1)^2$  and each class has $p$ elements.

\item \textit{$p\nmid i$ and $p\mid j$}. Suppose $j=tp$, where $1\leq t\leq p-1$. By Equation \ref{c4} and the fact that $a^{p^2}\in Z$, $(a^ub^v)^{-1}(a^{i}b^{tp})(a^ub^v)$ $=$ $a^{(utp-vi)p}a^{i}b^{tp}z_1=a^{-vip}a^{i}b^{tp}z_2$.  By Theorem \ref{thm3.5}, $|C_G(a^{i}b^{tp})|=p^2|Z|$ and so  $|(a^{i}b^{tp})^G|=p^2$. Now a similar argument as the Case 2,
\begin{equation}\label{c7}
(a^{i}b^{tp})^G=\{a^{kp}a^{i}b^{tp}z_{(k,l)}\mid 0\leq k,l\leq p-1,z_{(k,l)}\in Z\}.
\end{equation}

Hence, for fixed positive integers $i$ and $t$, there are $n$ distinct conjugacy classes as $(a^{i}b^{tp}z_{r_1})^G$, $\ldots$, $(a^{i}b^{tp}z_{r_n})^G$ in which $n=\frac{|Z|}{p}$.
Since $1\leq i\leq p^2-1$ and $p\nmid i$, there are integers $k'$ and $i'$ such that $i=k'p+i'$ and $1\leq i' \leq p-1$. By Equation \ref{c7},
\begin{eqnarray*}
(a^{i}b^{tp})^G & = & (a^{k'p+i'}b^{tp})^G \\
               & = & \{a^{kp}a^{k'p+i'}b^{tp}z_{(k,l)}\mid 0\leq k,l\leq p-1,z_{(k,l)}\in Z\} \\
               & = & \{a^{(k+k')p}a^{i'}b^{tp}z_{(k,l)}\mid 0\leq k,l\leq p-1,z_{(k,l)}\in Z\}.
\end{eqnarray*}
This shows that there exists $z_{(i',t)}\in Z$ such that $(a^{i}b^{tp})^G=(a^{i'}b^{tp}z_{(i',t)})^G$.  Since $1\leq i',t\leq p-1$,  the number of the distinct conjugacy classes is equal to $n(p-1)^2$ and each class has $p^2$ elements.

\item \textit{$p\mid i$ and $p\nmid j$}. Suppose $i=sp$ such that $1\leq s\leq p-1$. By Equation \ref{c4} and the fact that $a^{p^2}\in Z$, we have $(a^ub^v)^{-1}(a^{sp}b^{j})(a^ub^v)$ $=$ $a^{(uj-vsp)p}a^{sp}b^{j}z_1$ $=$ $a^{ujp}a^{sp}b^{j}z_2$ $=$ $a^{(uj+s)p}b^{j}z_2$. In this case, there is no new conjugacy class and the argument is similar to the Case 4.

\item \textit{$p\nmid i$ and $p\nmid j$}.   By Equation \ref{c4}, $(a^ub^v)^{-1}(a^{i}b^{j})(a^ub^v)=a^{(uj-vi)p}a^{i}b^{j}z$.
By Theorem \ref{thm3.5}, $|C_G(a^{i}b^{j})|=p^2|Z|$ and so $|(a^{i}b^{j})^G|=p^2$. Now by a similar argument as in the Case 2, we can see that
\begin{equation}\label{c8}
(a^{i}b^{j})^G=\{a^{kp}a^{i}b^{j}z_{(k,l)}\mid 0\leq k,l\leq p-1,z_{(k,l)}\in Z\}.
\end{equation}
Hence, for fixed positive integers $i$ and $j$, there are $n$ distinct conjugacy classes $(a^{i}b^{j}z_{r_1})^G,\cdots,(a^{i}b^{j}z_{r_n})^G$, where  $n=\frac{|Z|}{p}$.
Since $1\leq i\leq p^2-1$ and $p\nmid i$,  there are integers $k'$ and $i'$ such that $i=k'p+i'$ and $1\leq i'\leq p-1$. Therefore, by Equation \ref{c8},
\begin{eqnarray*}
(a^{i}b^{j})^G & = & (a^{k'p+i'}b^{j})^G \\
               & = & \{a^{kp}a^{k'p+i'}b^{j}z_{(k,l)}\mid 0\leq k,l\leq p-1,z_{(k,l)}\in Z\} \\
               & = & \{a^{(k+k')p}a^{i'}b^{j}z_{(k,l)}\mid 0\leq k,l\leq p-1,z_{(k,l)}\in Z\}.
\end{eqnarray*}
Thus, there exists $z_{(i',j)}\in Z$ such that $(a^{i}b^{j})^G=(a^{i'}b^{j}z_{(i',j)})^G$. Since $1\leq i'\leq p-1$, $1\leq j\leq p^2-1$ and $p\nmid j$, the number of distinct conjugacy classes  is equal to $np(p-1)^2$ and each class has $p^2$ elements.
\end{enumerate}
\end{enumerate}

\begin{center}
\begin{table}
\begin{tabular}{ |c|ll|c|l| }
\hline
 Type & \multicolumn{2}{|c|}{ The Representative of Conjugacy Classes }  & $|x^G|$ & $\#$ Conjugacy Classes \\ \hline
$1$ & $a^{sp}$ & $1\leq s\leq p-1$ & $p$  & $n(p-1)$ \\ \hline
$2$ & $a^i$ & $1\leq i\leq p-1$ & $p^2$  & $n(p-1)$ \\ \hline
$3$ & $b^{tp}$ & $1\leq t\leq p-1$ & $p$ & $n(p-1)$ \\ \hline
$4$ & $b^j$ & $1\leq j\leq p^2-1,j\overset{p}{\not\equiv} 0$ & $p^2$  & $np(p-1)$ \\ \hline
$5$ & $a^{sp}b^{tp}$ & $1\leq s,t\leq p-1$ & $p$  & $n(p-1)^2$ \\ \hline
$6$ & $a^ib^{tp}$ & $1\leq i,t\leq p-1$ & $p^2$  & $n(p-1)^2$ \\ \hline
$7$ & $a^{sp}b^j$ & $1\leq s\leq p-1,1\leq j\leq p^2-1,j\overset{p}{\not\equiv} 0$ & \multicolumn{2}{|c|}{The same of Type 4} \\ \hline
$8$ & $a^ib^j$ & $1\leq i\leq p-1,1\leq j\leq p^2-1,j\overset{p}{\not\equiv} 0$ & $p^2$ & $np(p-1)^2$ \\ \hline
\end{tabular}
\caption{Conjugacy Classes of $G$ such that $n=\frac{|Z|}{p}$.}\label{j3.2}
\end{table}
\end{center}

We now investigate the commuting conjugacy class graph of $G$. To do this, it is enough to determine the relationship between conjugacy classes in Table \ref{j3.2}. We consider the following cases:
\begin{enumerate}
\item Suppose $a^{s_1p}$ and $a^{s_2p}$ are the representatives of two classes of Type 1 in Table 2.
It is clear that $a^{s_1p}a^{s_2p}=a^{s_2p}a^{s_1p}$ and so all conjugacy classes of Type 1  are commuting together. Hence, the commuting conjugacy class graph has a subgraph isomorphic to $K_{n(p-1)}$. We now determine the relationship between conjugacy classes of this and another types. For this purpose, we suppose $a^{sp}$ is a representative of a conjugacy class of Type 1 and $a^ub^vz$ is an arbitrary element of $G$ such that they are commuted to each other. This means that $(a^{sp})(a^ub^vz) = (a^ub^vz)(a^{sp})$ and hence $a^{sp}b^v = b^va^{sp}$.
By Lemma \ref{lem3.2.5}(a), the last equality is true if and only if $p\mid v$ or $v=tp$ in which $0\leq t\leq p-1$. Then $a^{sp}$ commutes with all elements of the form $a^ub^{tp}z$. Therefore, the following cases can be occurred:
\begin{enumerate}
\item $t=0$. Then the conjugacy classes of Types 1 and 2 in Table 2 are commuting together.

\item \textit{$t\neq 0$ and $u=0$}. In this case, the conjugacy classes of Types 1 and 3 in Table 2 are commuting to each other.

\item \textit{$u,t\neq 0$ and $p\mid u$}. Then the conjugacy classes of Types 1 and 5 in Table 2 are commuting together.

\item \textit{$t\neq 0$ and $p\nmid u$}. By this condition, the conjugacy classes of Types 1  and 6  are commuting together.
\end{enumerate}

\item Suppose $a^{i_1}$ and $a^{i_2}$ are the representatives of two conjugacy classes of Type 2.
Obviously, $a^{i_1}a^{i_2}=a^{i_2}a^{i_1}$ and so all conjugacy classes of Type 2  are commuted to each other in the group. Hence, the commuting conjugacy class graph has a clique of size $n(p-1)$. We now determine the relation between this  and another types. To do this, we suppose that $a^i$ is a representative of a conjugacy class of Type 2  and $a^ub^vz$ is an arbitrary element of $G$ such that $(a^i)(a^ub^vz) = (a^ub^vz)(a^i)$. Thus,  $a^ib^v = b^va^i$.
Since $p\nmid i$, then By Lemma (\ref{lem3.2.5})(a) the last equality is satisfied if and only if $v=0$. Then $a^i$ commutes with all elements of the form $a^uz$ which shows that the conjugacy classes of Type 2  are commuting only with conjugacy classes of Type 1.

\item Suppose $b^{t_1p}$ and $b^{t_2p}$ are representatives of two conjugacy classes of Type 3.
It is clear that $b^{t_1p}b^{t_2p}=b^{t_2p}b^{t_1p}$. This shows that all conjugacy classes of Type 3 are commuting together and hence the commuting conjugacy class  graph has a subgraph isomorphic to $K_{n(p-1)}$. We now determine the relation between this and another types. Suppose $b^{tp}$ is a representatives of a conjugacy class of Type 3  and $a^ub^vz$ is an arbitrary element of $G$ such that $(b^{tp})(a^ub^vz) = (a^ub^vz)(b^{tp})$. Thus, $b^{tp}a^u = a^ub^{tp}$. By Lemma \ref{lem3.2.5}(a), the last equality is satisfied if and only if $p\mid u$ or $u=sp$ in which $0\leq s\leq p-1$. Then $b^{tp}$ commutes with all elements in the form of $a^{sp}b^vz$. Again, the following cases can be occurred:

\begin{enumerate}
\item \textit{$s=0$}. Then the conjugacy classes of Types 3 and 4 in Table 2 are commuting together.

\item \textit{$s\neq 0$ and $v=0$}. Under this condition, the conjugacy classes of Types 3 and 1 in Table 2 are commuting together.

\item \textit{$v,s\neq 0$ and $p\mid v$}. In this case, the conjugacy classes of Types 3 and 5 in Table 2 are commuting to each other.

\item \textit{$s\neq 0$ and $p\nmid v$}. The classes classes of Types 3 and 7 in Table 2 are commuting together.
\end{enumerate}

\item Suppose $b^{j_1}$ and $b^{j_2}$ are representatives of two conjugacy classes of Type 4.
Since $b^{j_1}b^{j_2}=b^{j_2}b^{j_1}$, all conjugacy classes of Type 4 are commuting together  and hence the commuting conjugacy class graph has a clique of size ${np(p-1)}$. We now determine the relationship between this and another types. To do this, we assume that $b^j$ is a representative of a conjugacy class of Type 4  and $a^ub^vz$ is an arbitrary element of $G$ such that $(b^j)(a^ub^vz) = (a^ub^vz)(b^j)$. Thus,  $b^ja^u = a^ub^j$.
Since $p\nmid j$, by Lemma \ref{lem3.2.5}(a) the last equality is satisfied if and only if $u=0$. Then $b^j$ commutes with all elements in the form of $b^vz$. So, the conjugacy classes of Type 4 are commuting only with the conjugacy classes of Type 3.

\item Suppose $a^{s_1p}b^{t_1p}$ and $a^{s_2p}b^{t_2p}$ are representatives of two conjugacy classes of Type 5.
By Equation \ref{e45}, $(a^{s_1p}b^{t_1p})(a^{s_2p}b^{t_1p})=(a^{s_2p}b^{t_2p})(a^{s_1p}b^{t_2p})$. This proves that all conjugacy classes of Type 5 are commuting to each other and hence the commuting conjugacy class graph has a clique of size ${n(p-1)^2}$. It is now enough to determine the relationship between this and conjugacy classes of types 6, 7 and 8. Suppose $a^{sp}b^{tp}$ and $a^ub^{vp}$ are representatives of conjugacy classes of Types 5 and 6, respectively, which are commuting together. By Lemma \ref{lem2.4}, $a^{sp}b^{tp}=(a^{ip}b^p)^{n_t}z_1$ and $a^ub^{vp}=(ab^{jp})^{n_u}z_2$. On the other hand, By Lemma \ref{lem3.2}, $C_G(a^{sp}b^{tp})=C_G(a^{ip}b^p)$ and $C_G(a^ub^{vp})=C_G(ab^{jp})$. We note that  $(a^{sp}b^{tp})(a^ub^{vp})=(a^ub^{vp})(a^{sp}b^{tp})$ if and only if $(a^{ip}b^p)(ab^{jp})=(ab^{jp})(a^{ip}b^p)$ if and only if $b^pa=ab^p$, contradicts by Lemma \ref{lem3.2.5}(a). Therefore, a conjugacy class of Type 5 is not adjacent with a conjugacy class of Type 6. We now assume that $a^{sp}b^{tp}$ and $a^ub^v$ are representatives of the conjugacy classes of Types 5 and 8, respectively, which are commute to each other. By Lemma \ref{lem2.4}, $a^{sp}b^{tp}=(a^{ip}b^p)^{n_t}z_1$ and $a^ub^v=(a^jb)^{n_v}z_2$. On the other hand, By Lemma \ref{lem3.2}, $C_G(a^{sp}b^{tp})=C_G(a^{ip}b^p)$ and $C_G(a^ub^v)=C_G(a^jb)$. By Equation 3.9, $(a^{sp}b^{tp})(a^ub^v)=(a^ub^v)(a^{sp}b^{tp})$ if and only if $(a^{ip}b^p)(a^jb)=(a^jb)(a^{ip}b^p)$ if and only if  $j \equiv i (mod \ p)$. Since $1\leq i\leq p-1$, the conjugacy classes of Types 5 and 8 can be divided into $p-1$ parts and hence  each part of Type 5 gives a complete subgraph of size ${n(p-1)}$.

\item Suppose  $a^{i_1}b^{t_1p}$ and $a^{i_2}b^{t_2p}$ are representatives of two commuting conjugacy classes of Type 6. By Lemma \ref{lem2.4}, $a^{i_1}b^{t_1p}$ $=$ $(ab^{s_1p})^{n_{i_1}}z_1$ and $a^{i_2}b^{t_2p}$ $=$ $(ab^{s_2p})^{n_{i_2}}z_2$. On the other hand, by Lemma \ref{lem3.2}, $C_G(a^{i_1}b^{t_1p})$ $=$ $C_G(ab^{s_1p})$ and $C_G(a^{i_2}b^{t_2p})$ $=$ $C_G(ab^{s_2p})$.
This proves that $(a^{i_1}b^{t_1p})(a^{i_2}b^{t_2p})$ $=$ $(a^{i_2}b^{t_2p})(a^{i_1}b^{t_1p})$ if and only if
$(ab^{s_1p})(ab^{s_2p})=(ab^{s_2p})(ab^{s_1p})$ if and only if $ab^{(s_2-s_1)p}=b^{(s_2-s_1)p}a$. By Lemma  \ref{lem3.2.5}(a), the last equality is satisfied if and only if $s_1=s_2$ and so in the commuting conjugacy class graph two conjugacy classes of Type 6 are adjacent, when their centralizers is equal. Also, by Theorem \ref{thm3.5}, the number of centralizers of Type 6 is $p-1$. So, the conjugacy classes of this type can be divided into $p-1$ parts. On the other hand, by Table \ref{j3.2},  the number of conjugacy classes of Type 6 are $n(p-1)^2$. Hence,
each part of Type 6 gives a complete subgraph isomorphic to $K_{n(p-1)}$.  It is now enough to find the relationship between conjugacy classes  of Types 6 and 8. To see this, we assume that  $a^ib^{tp}$ and $a^ub^v$ are representatives of  conjugacy classes of Types 6 and 8, respectively, such that they are commuting together. By Lemma \ref{lem2.4}, $a^ib^{tp}=(ab^{sp})^{n_i}z_1$ and $a^ub^v=(ab^j)^{n_u}z_2$. On the other hand, by Lemma \ref{lem3.2}, $C_G(a^ib^{tp})=C_G(ab^{sp})$ and $C_G(a^ub^v)=C_G(ab^j)$. Therefore, $(a^ib^{tp})(a^ub^v)=(a^ub^v)(a^ib^{tp})$ if and only if $(ab^{sp})(ab^j)=(ab^j)(ab^{sp})$ if and only if $b^{sp-j}a=ab^{sp-j}$. By Lemma \ref{lem3.2.5}(a), the last equality is satisfied if and only if $p^2\mid (sp-j)$ or $p\mid j$ which is impossible. This shows that in the commuting conjugacy class graph, the conjugacy classes of Type 6 is not adjacent with any conjugacy classes of Type 8.

\item Suppose  $a^{i_1}b^{j_1}$ and $a^{i_2}b^{j_2}$ are representatives of two commuting conjugacy classes of Type 8 such that they are commuting together. By Lemma \ref{lem2.4}, $a^{i_1}b^{j_1}=(a^sb)^{n_{j_1}}z_1$ and $a^{i_2}b^{j_2}=(a^tb)^{n_{j_2}}z_2$.
On the other hand, by Lemma \ref{lem3.2}, $C_G(a^{i_1}b^{j_1})=C_G(a^sb)$ and $C_G(a^{i_2}b^{j_2})=C_G(a^tb)$.
Therefore, $(a^{i_1}b^{j_1})(a^{i_2}b^{j_2})=(a^{i_2}b^{j_2})(a^{i_1}b^{j_1})$ if and only if $(a^sb)(a^tb)=(a^tb)(a^sb)$ if and only if $a^{s-t}b=ba^{s-t}$. By Lemma  \ref{lem3.2.5}(a), the last equality is satisfied if and only if $s=t$ and so in the commuting conjugacy class graph two conjugacy classes of Type 8 are adjacent, when their centralizers is equal. Also, by Theorem \ref{thm3.5}, the number of centralizers of Type 8 is equal to $p(p-1)$ and so the conjugacy classes of this type can be divided into $p(p-1)$ parts. By Table \ref{j3.1},  the number of conjugacy classes of Type 8 is  $np(p-1)^2$ and hence each part of Type 8 gives a complete subgraph of size ${n(p-1)}$. Therefore, we can see that in commuting conjugacy class graph of $G$, the conjugacy classes of Type 8  are adjacent only with conjugacy classes of Type 5. Finally, we know that the conjugacy classes of Types 5 and 8 can be divided into $p-1$ and $p(p-1)$ parts, respectively. Therefore, every conjugacy class of Type 5 is adjacent with $p$ conjugacy classes of Type 8.
\end{enumerate}

By above discussion, the commuting conjugacy class graph of $G$ is  a connected graph with $n(p-1)(p+1)^2$ vertices. Suppose $\mathcal{M}_2$ is the graph depicted in Figure \ref{fig4}. Then the commuting conjugacy class graph of $G$ can be written as a $\mathcal{M}_2$-join, i.e.
\[\Gamma(G)=\mathcal{M}_2[\overbrace{K_{n(p-1)},\cdots,K_{n(p-1)}}^{p^2+p+1},K_{np(p-1)}].\]
The  graph  $\Gamma(G)$ is depicted in Figure \ref{fig6}.
\begin{figure}[h]
\begin{picture}(400,300)
\scalebox{.95}{
\definecolor{qqqqff}{rgb}{0,0,0}
\begin{tikzpicture}[line cap=round,line join=round,>=triangle 45,x=1.2cm,y=1.2cm]
\clip(-5.66,-4.54) rectangle (15.68,4.94);
\draw [shift={(0.55,1.46)}] plot[domain=-0.37:3.87,variable=\t]({1*1.65*cos(\t r)+0*1.65*sin(\t r)},{0*1.65*cos(\t r)+1*1.65*sin(\t r)});
\draw [shift={(0.56,1.38)},dash pattern=on 2pt off 2pt]  plot[domain=3.84:5.96,variable=\t]({1*1.6*cos(\t r)+0*1.6*sin(\t r)},{0*1.6*cos(\t r)+1*1.6*sin(\t r)});
\draw [shift={(1.69,2.22)},dash pattern=on 2pt off 2pt]  plot[domain=0.11:0.71,variable=\t]({1*1.46*cos(\t r)+0*1.46*sin(\t r)},{0*1.46*cos(\t r)+1*1.46*sin(\t r)});
\draw [shift={(1.84,0.98)},dash pattern=on 2pt off 2pt]  plot[domain=5.51:6.1,variable=\t]({1*1.24*cos(\t r)+0*1.24*sin(\t r)},{0*1.24*cos(\t r)+1*1.24*sin(\t r)});
\draw [shift={(-0.42,0.53)},dash pattern=on 2pt off 2pt]  plot[domain=3.46:4.02,variable=\t]({1*1.67*cos(\t r)+0*1.67*sin(\t r)},{0*1.67*cos(\t r)+1*1.67*sin(\t r)});
\draw [shift={(-0.62,1.93)},dash pattern=on 2pt off 2pt]  plot[domain=2.39:2.93,variable=\t]({1*1.55*cos(\t r)+0*1.55*sin(\t r)},{0*1.55*cos(\t r)+1*1.55*sin(\t r)});
\draw (-1,2)-- (-2.22,1.76);
\draw (-2.14,2.26)-- (-1,2);
\draw (-1.75,2.98)-- (-1,2);
\draw (0.18,4.04)-- (0.37,3.09);
\draw (2.34,3.54)-- (1.9,2.38);
\draw (2.82,3.2)-- (1.9,2.38);
\draw (3.14,2.38)-- (1.9,2.38);
\draw (3.14,1.22)-- (2.09,0.86);
\draw (3.05,0.76)-- (2.09,0.86);
\draw (2.73,0.12)-- (2.09,0.86);
\draw (-1,-1)-- (-0.67,0.36);
\draw (-1.47,-0.74)-- (-0.67,0.36);
\draw (-2,0)-- (-0.67,0.36);
\draw (1.9,2.38)-- (-1,2);
\draw (-1,2)-- (2.09,0.86);
\draw (2.09,0.86)-- (-0.67,0.36);
\draw (1.9,2.38)-- (-0.67,0.36);
\draw (0.37,3.09)-- (-0.67,0.36);
\draw (0.37,3.09)-- (2.09,0.86);
\begin{scriptsize}
\fill [color=qqqqff] (-1,2) circle (1.5pt);
\draw[color=qqqqff] (-0.9,1.8) node {$1$};
\fill [color=qqqqff] (-2.22,1.76) circle (1.5pt);
\draw[color=qqqqff] (-2.4,1.7) node {$2$};
\fill [color=qqqqff] (-2.14,2.26) circle (1.5pt);
\draw[color=qqqqff] (-2.3,2.2) node {$6_1$};
\fill [color=qqqqff] (-1.75,2.98) circle (1.5pt);
\draw[color=qqqqff] (-1.7,3.17) node {$6_{p-1}$};
\fill [color=qqqqff] (0.37,3.09) circle (1.5pt);
\draw[color=qqqqff] (0.4,2.8) node {$3$};
\fill [color=qqqqff] (0.18,4.04) circle (1.5pt);
\draw[color=qqqqff] (0.2,4.22) node {$4$};
\fill [color=qqqqff] (1.9,2.38) circle (1.5pt);
\draw[color=qqqqff] (1.9,2.2) node {$5_1$};
\fill [color=qqqqff] (2.34,3.54) circle (1.5pt);
\draw[color=qqqqff] (2.45,3.73) node {$8_{1_1}$};
\fill [color=qqqqff] (2.82,3.2) circle (1.5pt);
\draw[color=qqqqff] (3.02,3.38) node {$8_{1_2}$};
\fill [color=qqqqff] (3.14,2.38) circle (1.5pt);
\draw[color=qqqqff] (3.4,2.3) node {$8_{1_p}$};
\fill [color=qqqqff] (2.09,0.86) circle (1.5pt);
\draw[color=qqqqff] (2.15,.6) node {$5_2$};
\fill [color=qqqqff] (3.14,1.22) circle (1.5pt);
\draw[color=qqqqff] (3.35,1.2) node {$8_{2_1}$};
\fill [color=qqqqff] (3.05,0.76) circle (1.5pt);
\draw[color=qqqqff] (3.3,0.7) node {$8_{2_2}$};
\fill [color=qqqqff] (2.73,0.12) circle (1.5pt);
\draw[color=qqqqff] (3.,0.0) node {$8_{2_p}$};
\fill [color=qqqqff] (-0.67,0.36) circle (1.5pt);
\draw[color=qqqqff] (-0.2,0.25) node {$5_{p-1}$};
\fill [color=qqqqff] (-1,-1) circle (1.5pt);
\draw[color=qqqqff] (-0.6,-1.2) node {$8_{(p-1)_1}$};
\fill [color=qqqqff] (-1.47,-0.74) circle (1.5pt);
\draw[color=qqqqff] (-1.6,-0.9) node {$8_{(p-1)_2}$};
\fill [color=qqqqff] (-2,0) circle (1.5pt);
\draw[color=qqqqff] (-2.,0.3) node {$8_{(p-1)_p}$};
\end{scriptsize}
\end{tikzpicture}
}
\end{picture}
\vspace*{-3cm}
\caption{The Graph $\Gamma(G)$, when $\frac{G}{Z(G)}$ is non-abelian and $\frac{G}{Z}\cong \mathbb{Z}_{p^2}\rtimes \mathbb{Z}_{p^2}$ .}\label{fig6}
\end{figure}
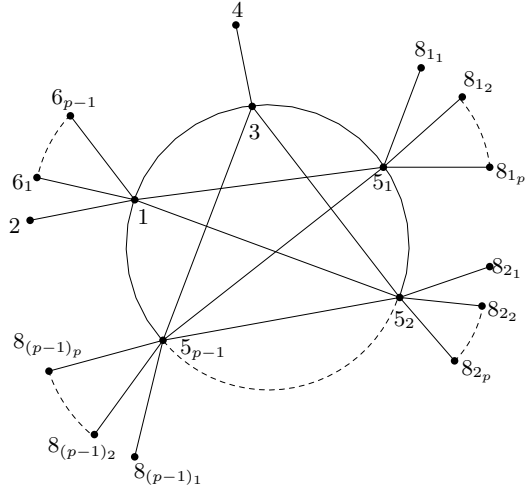
\end{proof}


\vskip 3mm

\noindent{\bf Acknowledgement.}  The research of the first author is partially supported by the University of Kashan under grant no. 364988/64.

\end{document}